\documentclass[final,hidelinks,onefignum,onetabnum]{siamart220329}

\usepackage{lipsum}
\usepackage{amsfonts}
\usepackage{graphicx}
\usepackage{epstopdf}
\usepackage{algorithmic}

\usepackage{subfig}
\usepackage{tikz}
\usetikzlibrary{calc,patterns,angles,quotes}

\ifpdf
  \DeclareGraphicsExtensions{.eps,.pdf,.png,.jpg}
\else
  \DeclareGraphicsExtensions{.eps}
\fi

\newsiamremark{remark}{Remark}
\newsiamremark{hypothesis}{Hypothesis}
\crefname{hypothesis}{Hypothesis}{Hypotheses}
\newsiamthm{claim}{Claim}

\headers{Shape-Newton Method for FBP Subject to Bernoulli Condition}{Y. Fan, J. Billingham, and K. van der Zee}

\title{A Shape-Newton Method for Free-boundary Problems Subject to The Bernoulli Boundary Condition}

\author{Yiyun Fan
\and John Billingham\thanks{School of Mathematical Sciences, University of Nottingham(\email{john.billingham@nottingham.ac.uk}, \email{KG.vanderZee@nottingham.ac.uk}).}
\and Kristoffer van der Zee\footnotemark[1]}

\usepackage{amsopn}

\ifpdf
\hypersetup{
  pdftitle={A Shape-Newton Method for Free-boundary Problems Subject to The Bernoulli Boundary Condition},
  pdfauthor={Yiyun Fan, John Billingham, and Kristoffer van der Zee}
}
\fi

\begin{document}

\maketitle

\begin{abstract}
We develop a shape-Newton method for solving generic free-boundary problems where one of the free-boundary conditions is governed by the nonlinear Bernoulli equation. The method is a Newton-like scheme that employs shape derivatives of the governing equations. In particular, we derive the shape derivative of the Bernoulli equation, which turns out to depend on the curvature in a nontrivial manner. The resulting shape-Newton method allows one to update the position of the free boundary by solving a special linear boundary-value problem at each iteration. We prove solvability of the linearised problem under certain conditions of the data. We verify the effectiveness of the shape-Newton approach applied to free-surface flow over a submerged triangular obstacle using a finite element method on a deforming mesh. We observe superlinear convergence behaviour for our shape-Newton method as opposed to the unfavourable linear rate of traditional methods.
\end{abstract}



\section{Introduction}
Free boundary problems have many applications in fluid mechanics, such as open-channel flow, fluid/solid interaction and hydrodynamics. Solving such problems is difficult, because the geometry of the domain needs to be determined together with other variables in this problem. A simplified but important model problem is the Bernoulli free-boundary problem, which considers a (linear) Dirichlet boundary condition, as well as a Neumann boundary condition on the free boundary~\cite{crank1987free, rumpf1997bernoulli}. This problem is not to be confused with the Bernoulli \emph{equation}, which is the pressure boundary condition in irrotational fluid mechanics, and which we will study in this paper. The nonlinearity of the Bernoulli equation poses an additional challenge to numerical algorithms.

There are several computational approaches to solving free-boundary problems. The first is to iteratively solve the boundary value problem with a single free-boundary condition for the field variables on a fixed approximated domain, and then update the free surface derived from the remaining free boundary condition (which was not included in the boundary value problem). These fixed-point type methods are called \emph{trial} methods, which converge linearly and cannot always find a solution. Details can be found, for example, in \cite{rumpf1997bernoulli, bouchon2005numerical,kuster2007fast}.

The second approach is to formulate a shape optimization problem to improve the convergence rate. This method aims to construct a boundary-value problem as the state problem with one free-boundary condition and formulate a cost function with the remaining free-boundary condition. This approach may require gradient information. The formulation and application of shape optimization to free boundary problems can be found in, e.g.~\cite{eppler2006efficient, haslinger1992optimal, haslinger2003shape, tiihonen1997shape,toivanen2008shape,van2003numerical}.

The third approach requires \emph{linearising} the whole system and applying a Newton-type method. The use of shape calculus and a Newton-type method is called the shape-Newton method. One linearisation method, called \emph{domain-map} linearisation, requires to transform the free-boundary problem to an equivalent boundary value problem on a fixed domain and then linearise the transformed problem with respect to the domain map~\cite{mejak1994numerical, VANDERZEE20112738}. An alternative way to linearise the free-boundary problem is to apply \emph{shape linearisation}~\cite{Shapes, sokolowski1992introduction}. K{\"a}rkk{\"a}inen and Tiihonen used this technique to solve Bernoulli free-boundary problems~\cite{karkkainen2004shape, karkkainen1999free}. The application to a more general Bernoulli free-boundary problem has been investigated in Van der Zee et al~\cite{van2013shape} by considering the whole problem in one weak form, and using $C^1$-continuous $B$-splines to represent discrete free boundaries, in order to allow the exact computation of the curvature in the shape derivatives. Montardini et al.~\cite{montardini2021isogeometric} extend this method by incorporating a collocation approach to update the boundary, and compare both methods by imposing Dirichlet or periodic boundary conditions on the vertical fixed boundary of the domain. The results show that collocation scheme has slightly worse accuracy but higher efficiency. 

In the current work, we derive the shape-Newton method for a free-boundary problem involving the nonlinear Bernoulli boundary condition on the free boundary. We use our approach to also re-derive the shape-Newton method for the simpler Bernoulli free-boundary problem (containing a Dirichlet boundary condition), which was obtained in~\cite{van2013shape} using a slightly different derivation.\footnote{We note that there is a typo for the strong form of the linearised problem in~\cite{van2013shape}. This mistake is rectified in this paper; see equations~\eqref{SF_laplace_D}--\eqref{CS_D}.} Similar to K{\"a}rkk{\"a}inen and Tiihonen, we set up two weak statements: One derived from the boundary value problem with the Neumann boundary condition, and the other from the remaining free boundary condition (Dirichlet condition or nonlinear Bernoulli condition).

A key result in our work is the 
shape derivative of the Bernoulli equation. 
It turns out that it has various equivalent expressions that are  surprisingly elegant: The primary result involves the normal derivative of the velocity squared ($|\nabla^2\phi|$), and we show in detail how this can be equivalently computed using only the tangential components of the velocity, suitably weighted by curvatures; see Section \ref{linearisation_bernoulli}.

We present our shape-Newton scheme in both strong and weak form, and without reference to any particular underlying discretisation. We study the  solvability of the linearised system in the continuous setting, that is, we establish coercivity of a suitable bilinear form under certain conditions of the data. We are also able to establish discrete solvability for a particular finite element approximation using deforming meshes, under certain conditions. We show numerical experiments involving open channel flow over a submerged triangle. We observe that the shape-Newton method converges superlinearly, and the results agree well with exact solutions and results from~\cite{dias1989open}.

The contents of this paper are arranged as follows. We first introduce the model problems either with the Dirichlet boundary condition or the Bernoulli equation on the free boundary in Section~\ref{chap:3.2}. In Section \ref{chap:3.3}, we derive the weak form for both problems. Then, we introduce some basic concepts about shape derivatives in Section~\ref{chap:3.4}. We carry out shape linearisation by applying Hadamard shape derivatives for the free-boundary problem in Section~\ref{chap:3.5}. In this Section, we also present the various equivalent expressions for the shape derivative of the Bernoulli equation. In Section~\ref{chap:3.6}, we present the Newton-like schemes, and present solvability results for the involved linearised systems (details in Appendix~\ref{app:solvability_dirichlet} and~\ref{app:solvability_bernoulli}). The finite element scheme using deforming meshes is given in Section~\ref{chap:3.8} (details of its discrete solvability in Appendix~\ref{app:solvability_discrete}), as well as numerical experiments. These are followed by Conclusions in Section~\ref{chap:3_conclusion}.

\section{Free-boundary Problem with Bernoulli or Dirichlet free-boundary condition}\label{chap:3.2}

We investigate the free boundary problem with either the Bernoulli condition or the Dirichlet condition on the free boundary. The Bernoulli condition is commonly used when considering steady, incompressible, and inviscid flow, but it is nonlinear, making the free-boundary problem more challenging to solve. To be general, the boundary conditions on the fixed boundaries are Robin boundary conditions.

\subsection{Free-boundary Problem With Bernoulli Condition}\label{chap3:bernoulli}
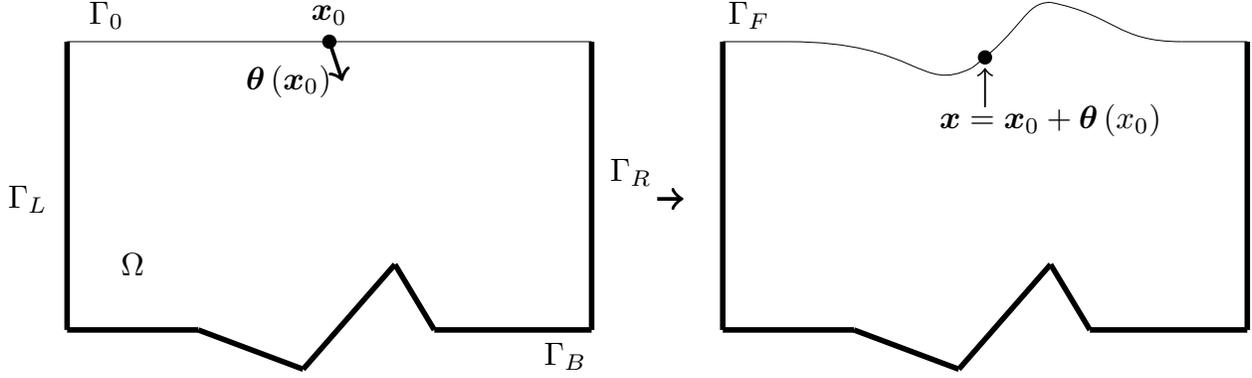
\begin{figure}[t!]
\centering
\resizebox{\columnwidth}{!}{%
\begin{tikzpicture}
\draw [very thick] (-8, 0) -- (-7,0);
\draw [very thick] (-7,0) -- (-6.2,-0.3);
\draw [very thick] (-5.5,0.5) -- (-6.2,-0.3);
\draw [very thick] (-5.5,0.5) -- (-5.2,0);
\draw [very thick] (-4,0) -- (-5.2,0);

\draw [very thick] (-8, 0) -- (-8,2.2);
\draw [very thick] (-4,0) -- (-4,2.2);

\draw [very thin] (-8, 2.2) -- (-4, 2.2);

\draw [->, thick] (-3.5,1) -- (-3.3,1);

\draw [very thick] (-3, 0) -- (-2,0);
\draw [very thick] (-2,0) -- (-1.2,-0.3);
\draw [very thick] (-0.5,0.5) -- (-1.2,-0.3);
\draw [very thick] (-0.5,0.5) -- (-0.2,0);
\draw [very thick] (1,0) -- (-0.2,0);

\draw [very thick] (-3, 0) -- (-3,2.2);
\draw [very thick] (1,0) -- (1,2.2);
\draw [very thin] (-3, 2.2) -- (-2.5, 2.2);
\draw [very thin] (0.5, 2.2) -- (1, 2.2);

\draw[very thin] (-2.5,2.2) .. controls (-1.5,2.2) and (-1.5,1.8) .. (-1.1,2);
\draw[very thin] (-0.5,2.5) .. controls (-0.7,2.5) and (-0.7,2.3) .. (-1.1,2);
\draw[very thin] (0.5,2.2) .. controls (0,2.2) and (0,2.4) .. (-0.5,2.5);

\node at (-2.8,2.4)[ inner sep=0pt, scale=0.7] {$\Gamma_F$};
\node at (-3.7,1.2)[ inner sep=0pt, scale=0.7] {$\Gamma_R$};
\node at (-8.3,1)[ inner sep=0pt, scale=0.7] {$\Gamma_L$};
\node at (-7.5,0.5)[ inner sep=0pt, scale=0.7] {$\Omega_0$};
\node at (-2.5,0.5)[ inner sep=0pt, scale=0.7] {$\Omega$};

\node at (-4.2,-0.2)[ inner sep=0pt, scale=0.7] {$\Gamma_B$};
\node at (-7.7,2.4)[ inner sep=0pt, scale=0.7] {$\Gamma_0$};

\draw [black,fill=black](-6.5,2.2) circle (.3ex); 
\draw [black,fill=black](-1.3,1.95) circle (.3ex); 
\node at (-7, 2.4)[ inner sep=0pt, scale=0.7] {$\boldsymbol x_0$};
\draw [->, thick] (-6.5,2.2) -- (-6.2, 1.95 )  node[left,scale=0.7,very thin] {$\boldsymbol\theta\left(\boldsymbol x_0\right)$};
\node at (-1.5, 1.6)[ inner sep=0pt, scale=0.7] {$\boldsymbol x = \boldsymbol x_0 + \boldsymbol\theta\left(x_0\right)$};
\draw [->, thin] (-1.5, 1.7) -- (-1.3,1.9);

\end{tikzpicture}%
}
\caption{The sketch of the parametrization of the free boundary $\Gamma_F$ by the displacement $\boldsymbol\theta\left(\boldsymbol x_0\right)$ with respect to the reference boundary $\Gamma_0$.}
\label{chap3:Illustration}
\end{figure}

The \break free-boundary problem with a Bernoulli condition can be stated as seeking an unknown domain $\Omega\subset\mathbb{R}^N$ ($N=2, 3$ or even $>3$), and a corresponding scalar potential function $\phi: \Omega\to\mathbb{R}$. For fluid problems, $\nabla\phi$ is then the velocity vector. The boundary $\partial\Omega$ contains a free boundary $\Gamma_F$, and the remainder $\partial\Omega\setminus\Gamma_F$, for example in the two-dimensional open-channel flow case, $\partial\Omega\setminus\Gamma_F$ contains a left boundary $\Gamma_L$ for inflow, a right boundary $\Gamma_R$ for outflow, and the bed $\Gamma_B$ which can have any reasonable shape. Figure \ref{chap3:Illustration} is an example of the domain and the parametrization of the free boundary $\Gamma_F$.

The problem can be written as
\begin{subequations}
\begin{eqnarray}
-\Delta\phi &=& f, \quad \text{ in } \Omega,\label{FB_L}\\
\partial_{\boldsymbol n}\phi &=& 0, \quad \text{ on }\Gamma_F,\label{FB_N}\\
a\left|\nabla\phi\right|^2+b x_N+c &=& 0, \quad \text{ on } \Gamma_F,\label{FB_BC}\\
\partial_{\boldsymbol n}\phi + \omega\phi &=& g + \omega h, \quad \text{ on }\partial\Omega\setminus\Gamma_F,\label{FB_Fixed}
\end{eqnarray} 
\end{subequations}
where $\partial_{\mathbf{n}}\left(\cdot\right) = \mathbf{n}\cdot\nabla\left(\cdot\right)$ is the normal derivative with $\mathbf{n}$ being the unit normal vector to the boundary pointing out the domain, and $x_N$ is the $N$-th component (vertical component) of vector~$\mathbf{x}$. The condition \eqref{FB_L} is the PDE for potential $\phi$, where $f: \mathbb{R}^N\to\mathbb{R}$ is a sufficiently smooth given function. The condition \eqref{FB_N} represents the kinematic condition on the free boundary. The condition \eqref{FB_BC} with real-valued constants $a, b,$ and $c$ represents the Bernoulli condition.%
\footnote{Because of~\eqref{FB_N}, the Bernoulli condition~\eqref{FB_BC} is a condition on~$|\nabla_\Gamma \phi|$, hence it can be thought of as a surface-eikonal equation~\cite{spira2004efficient, fu2011fast}.} 
In the standard case, $a=\frac{1}{2}$, $b$ is the gravitational acceleration, and $c = \frac{p_{\infty}}{\rho_0}$, where $p_{\infty}$ is the external pressure and $\rho_0$ is the constant density of the fluid.

We consider general Robin boundary conditions \eqref{FB_Fixed} on $\partial\Omega\setminus\Gamma_F$ where $\omega\geq0$, and $g, h:\mathbb{R}^N\to\mathbb{R}$ are sufficiently smooth given functions. Thus we can approximate either a Neumann or Dirichlet-type condition depending on the value of $\omega$: the Neumann boundary condition, obtained when $\omega=0$, usually represents the kinematic condition, where the perpendicular fluid velocity is zero on the free or solid boundary. On the other hand, choosing $\omega\to\infty$ yields the Dirichlet boundary condition $\phi = h$. Furthermore, it is possible to impose mixed boundary conditions by choosing various values of $\omega$ on different parts of the boundaries (e.g. $\Gamma_L$, $\Gamma_R$ and $\Gamma_B$). 

We assume that for suitable data $f$, $g$, $h$ and $a, b, c, \omega$, there is a nontrivial and sufficiently-smooth solution pair $(\Gamma_F, \phi)$. The wellposedness of a free-boundary problem is studied in, for example, \cite{hamilton1982inverse, saavedra1991variational}. By introducing a vector field $\boldsymbol\theta: \Gamma_0\to\mathbb{R}^N$, the displacement of the free boundary with respect to the referenced boundary $\Gamma_0$ (of constant height $\bar{x}_N$) can be defined as
\begin{eqnarray}
\Gamma_F :=\left\{\mathbf{x}\in\mathbb{R}^N|\mathbf{x} = \mathbf{x}_0 + \boldsymbol\theta\left(\mathbf{x}_0\right), \forall \mathbf{x}_0\in\Gamma_0\right\} 
\end{eqnarray}
to parametrize the domain $\Omega$ and the free boundary $\Gamma_F$, as shown in Figure  \ref{chap3:Illustration}. This allows us to think of the problem \eqref{FB_L}--\eqref{FB_Fixed} in terms of the solution pair $\left(\boldsymbol\theta, \phi\right)$. 

\begin{remark}[Fixed free boundary at inflow]
    The part of the free boundary $\Gamma_F$ corresponding to inflow is assumed to be fixed in our work. That means in the two-dimensional case that the left node $\boldsymbol{x}_L$ on the free boundary $\Gamma_F$ is fixed, i.e. $\boldsymbol{\theta}(\boldsymbol{x_L})=\boldsymbol0$.
    \label{remark_compatibility1}
\end{remark}
\begin{remark}[Data compatibility I] 
    There is a compatibility requirement on the parameters defining the Bernoulli and Robin boundary conditions, \eqref{FB_BC} and \eqref{FB_Fixed}. It is well-known that the Bernoulli condition prescribes the conservation of energy along streamlines, therefore the total energy prescribed by the Bernoulli condition should match the value of the total energy given by the Robin boundary condition at the corresponding upstream (and downstream) coordinates. For example, in the two-dimensional case with the geometry corresponding to the illustration in Fig.\ref{chap3:Illustration}, the choice $\omega|_{\boldsymbol{x}_L}=0$, $g|_{\boldsymbol{x}_L}=-\cos\alpha$, and $a=\frac{1}{2}F^2,bx_N=1,c=-\frac{1}{2}F^2-1$ are compatible for inflow angle~$\alpha$ and Froude number~$F$. 
    \label{remark_compatibility2}
\end{remark}
\begin{remark}[Data compatibility II]
    Further to the compatibility requirements in the previous remark, we will also require the following, which will ensure the linearised operator is well-posed (see Remark \ref{remark_solve_dirichlet}):
    
    Let $\boldsymbol{\tau}$ denote the unit vector at $\partial\Gamma_F$, tangent to and outward from $\Gamma_F$, and normal to $\partial\Gamma_F$. We require $\phi$ to have $\nabla_{\Gamma}\phi\cdot\boldsymbol{\tau} = 0$ on those parts of $\partial\Gamma_F$ that do not touch a part of $\partial\Omega\setminus\Gamma_{F}$ where a Dirichlet boundary condition is imposed, or where $\boldsymbol{\theta}$ is fixed. An example situation for which this is satisfied is a flow in an open channel with non-homogeneous Neumann boundary condition on the lateral boundary and Dirichlet boundary condition on the outflow. An alternative situation is where the Dirichlet boundary condition holds on all over $\partial\Omega$.
    \label{remark_compatibility3}
\end{remark}

\subsection{Free-boundary Problem with Dirichlet Boundary Condition} \label{chap3:dirichlet}

A more simple model problem is introduced by replacing the Bernoulli condition with the Dirichlet condition on the free boundary. The dependence on $\phi$ is now linear:
\begin{subequations}
\begin{eqnarray}
-\Delta\phi &=& f, \quad \text{ in } \Omega,\label{FBD_L}\\
\partial_{\boldsymbol n}\phi &=& 0, \quad \text{ on }\Gamma_F,\label{FBD_N}\\
\phi &=& h, \quad \text{ on } \Gamma_F,\label{FBD_D}\\
\partial_{\boldsymbol n}\phi + \omega\phi &=& g + \omega h, \quad \text{ on }\partial\Omega\setminus\Gamma_F,\label{FBD_LBC}
\end{eqnarray}  
\end{subequations}
where $g$ and $h$ are assumed to be sufficiently smooth on $\mathbb{R}^N$ (e.g. $g\leq0$ on $\Gamma_L$ makes sense for the inflow). 

By choosing $\omega\to\infty$ on $\partial\Omega\setminus\Gamma_F$ (i.e., Dirichlet boundary condition), this problem becomes the classical
problem for an ideal fluid, called the Bernoulli free-boundary problem \cite{rumpf1997bernoulli}. 

\section{The Weak Form}\label{chap:3.3}

We will first find weak forms of both free-boundary problems in order to apply shape-calculus techniques to linearise these problems, and  subsequently propose Newton-like schemes. To allow shape linearisation, the test functions will be more regular than what is usually assumed. Hence, let $v\in V :=H^2(\mathbb{R}^N)$ and $w\in W := {H}^2(\mathbb{R}^N)$ be sufficiently smooth test functions.

Since the only difference between the two free-boundary problems in Section \ref{chap:3.2} is the Bernoulli condition and the Dirichlet condition on the free boundary, the first weak form in the domain $\Omega$ is the same in both situations. It can be obtained by multiplying the Laplacian equation (\eqref{FB_L} or \eqref{FBD_L}) by the test function $v\in V$ and integrating over $\Omega$, then applying the Green's formula with the Robin boundary conditions on $\partial\Omega\setminus\Gamma_F$ and (homogeneous) Neumann boundary condition on $\Gamma_F$, yielding
\begin{equation}
\mathcal{R}_1\left(\left(\boldsymbol\theta,\phi\right); v\right) = 0, \quad \forall v\in V,\label{nonlinear_1}
\end{equation}
where the semilinear form $\mathcal{R}_1\left(\left(\boldsymbol\theta,\phi\right); v\right)$ is defined as
\begin{eqnarray}
\mathcal{R}_1\left(\left(\boldsymbol\theta,\phi\right); v\right) &=& 
\int_{\Omega}\nabla\phi\cdot\nabla v  d \Omega -\int_{\partial\Omega\setminus\Gamma_F}\left(g+\omega h-\omega\phi\right)v\text{d}\Gamma - \int_{\Omega}fv d \Omega\label{weak form 1}.
\end{eqnarray}

The second weak form can be derived by multiplying the remaining free-boundary condition by the test function $w\in W$ and integrating over $\Gamma_F$,
\begin{equation}
\mathcal{R}_2\left(\left(\boldsymbol\theta, \phi\right); w\right) = 0, \quad \forall w\in W,\label{nonlinear_2}
\end{equation}
with the definition of the semilinear form $\mathcal{R}_2\left(\left(\boldsymbol\theta, \phi\right); w\right)$ as
\begin{equation}
\mathcal{R}_2\left(\left(\boldsymbol\theta, \phi\right); w\right) = \int_{\Gamma_F}\left(\text{B.C}\right)w  d \Gamma, \label{weak form 2}
\end{equation}
where $\left(\text{B.C}\right)$ can either be the left hand side of Bernoulli condition \eqref{FB_BC} or $(\phi-h)$ in case of the Dirichlet condition \eqref{FBD_D}.

Given some approximation $(\hat{\boldsymbol{\theta}}, \hat{\phi})$, the exact Newton method for an update \hfill\break$\left(\delta\boldsymbol\theta, \delta\phi\right)$, in weak form, would be 
\begin{subequations}
\begin{eqnarray}
\left<\partial_{(\boldsymbol\theta, \phi)}\mathcal{R}_1\left(\left(\hat{\boldsymbol\theta}, \hat{\phi}\right);v\right), \left(\delta\boldsymbol\theta, \delta\phi\right)\right> &=& - \mathcal{R}_1\left(\left(\hat{\boldsymbol\theta}, \hat{\phi}\right); v\right) \quad \forall v\in V,\nonumber \\
\left<\partial_{(\boldsymbol\theta, \phi)}\mathcal{R}_2\left(\left(\hat{\boldsymbol\theta}, \hat{\phi}\right);w\right), \left(\delta\boldsymbol\theta, \delta\phi\right)\right> &=& - \mathcal{R}_2\left(\left(\hat{\boldsymbol\theta}, \hat{\phi}\right); v\right) \quad \forall w\in W.\nonumber
\end{eqnarray}
\end{subequations}
We now study the shape derivatives, which are present in the above left-hand side.

\section{Shape Derivatives}\label{chap:3.4}

The linearisation of $\mathcal{R}_1\left(\left(\boldsymbol\theta,\phi\right); v\right)$ and $\mathcal{R}_2\left(\left(\boldsymbol\theta, \phi\right); w\right)$ needs the differentiation of the weak forms with respect to the geometry, where the geometry itself is treated as a variable. Thus the shape derivatives are applied to a given domain, which requires some appropriate smoothness assumptions. 

The weak forms \eqref{weak form 1} and \eqref{weak form 2} contain domain integrals $\int_{\Omega}\left(\cdot\right) d \Omega$ and boundary integrals $\int_{\Gamma_F}\left(\cdot\right) d \Gamma$. The shape derivatives for a domain integral and a boundary integral can be obtained by the Hadamard formulas \cite{Shapes,  sokolowski1992introduction}:

\begin{theorem}[Shape derivative of domain integral]
Suppose $\phi\in W^{1,1}\left(\mathbb{R}^N\right)$, where 
$$W^{1,1}\left(\mathbb{R}^N\right) = \{f\in L^1(\mathbb{R}^N): \nabla f\in L^1(\mathbb{R}^N)^N\},$$
and $\Omega$ is an open and bounded domain with boundary $\Gamma = \partial \Omega$ of class $C^{0,1}$. Consider the domain integral
\begin{equation}
J\left(\Omega\right) = \int_{\Omega}\phi  d \Omega.\nonumber
\end{equation}
Then its shape derivative with respect to the perturbation $ \delta\boldsymbol \theta\in C^{0,1}\left(\mathbb{R}^N;\mathbb{R}^N\right)$ is given by%
\footnotemark{}%
\begin{equation}
\left< d J\left(\Omega\right), \delta\boldsymbol\theta\right> = \int_{\Gamma}\phi\delta \boldsymbol\theta\cdot\boldsymbol n d \Gamma,\nonumber
\end{equation}
where $\boldsymbol n$ denotes the outward normal derivative to $\Omega$.
\label{Shape_domain}
\end{theorem}
\footnotetext{\label{footnote_linearise}$\left< d J\left(\Omega\right), \delta\boldsymbol\theta\right> := \lim_{t\rightarrow 0} \frac{1}{t}\big(J(\Omega + \delta\boldsymbol\theta) - J(\Omega)\big)$, see e.g.,~\cite[Chapter~9]{Shapes}. In practise, $\delta \boldsymbol{\theta}$ is only needed on $\Gamma = \partial \Omega$, instead of on the whole~$\mathbb{R}^{N}$. Any extension of~$\delta \boldsymbol{\theta}\in C^{0,1}(\Gamma;\mathbb{R}^N)$ into $C^{0,1}(\mathbb{R}^N;\mathbb{R}^N)$ would suffice since $\left< d J\left(\Omega\right), \delta\boldsymbol\theta\right>$ does not depend on the particular extension used.}

\begin{theorem}[Shape derivative of boundary integral]
Suppose $\phi\in W^{2,1}\left(\mathbb{R}^N\right)$, where
$$W^{2,1}\left(\mathbb{R}^N\right) = \{f\in L^1(\mathbb{R}^N):D^kf\in L^1(\mathbb{R}^N) \text{ for any multi-index~$k$ with } |k|\leq 2\},$$
and $\Omega$ is an open and bounded domain with boundary $\Gamma = \partial \Omega$ of class $C^{1,1}$. Consider the boundary integral
\begin{equation}
J\left(\Omega\right)=\int_{\Gamma}\phi  d \Gamma.\nonumber
\end{equation}
Then its shape derivative with respect to the perturbation $ \delta\boldsymbol\theta\in C^{0,1}\left(\mathbb{R}^N;\mathbb{R}^N\right)$ is given by\footref{footnote_linearise}
\begin{equation}
\left< d J\left(\Omega\right), \delta\boldsymbol\theta\right>  = \int_{\Gamma}\left(\partial_{\boldsymbol n}\phi+\kappa\phi\right)\delta\boldsymbol\theta\cdot\boldsymbol n d \Gamma,\nonumber
\end{equation}
where $\boldsymbol n$ denotes the normal vector to $\Gamma$ and $\kappa$ is the (additive) curvature of $\Gamma$.
\label{Shape_boundary_closed}
\end{theorem}

\begin{remark}
    The shape derivative of boundary integral for the open boundary (see \cite[Eq. (5.48)]{walker2015shapes})  is: \begin{equation}
        \left< d J\left(\Omega\right), \delta\boldsymbol\theta\right>  = \int_{\Gamma}\left(\partial_{\boldsymbol n}\phi+\kappa\phi\right)\delta\boldsymbol\theta\cdot\boldsymbol n d \Gamma + \int_{\partial\Gamma}\phi\boldsymbol{\tau}\cdot\delta\boldsymbol{\theta}ds,\nonumber\label{Shape_boundary_open}
    \end{equation}
    where $\boldsymbol{\tau}$ is defined in Remark \ref{remark_compatibility3}.
    
\end{remark}
\begin{remark}[Piecewise-smooth free boundary]
    When~$\Gamma$ is piecewise smooth, additional jump terms should be included in the boundary-integral shape derivative; see e.g.,~\cite[Ch.~3.8]{sokolowski1992introduction}. \label{remark_piecewise}
\end{remark}

\section{Linearisation}\label{chap:3.5}
The linearisation of $\mathcal{R}_1\left(\left(\boldsymbol\theta,u\right);v\right)$ and $\mathcal{R}_2\left(\left(\boldsymbol\theta,\phi\right);w\right)$ at an approximation pair $\left(\hat{\boldsymbol\theta}, \hat{\phi}\right)$ close to the exact solution~$(\boldsymbol{\theta}^*, \phi^*)$ can be derived from the partial derivative of the weak forms with respect to $\phi$ and $\boldsymbol\theta$. We proceed formally when obtaining our linearisation: We assume that $\hat{\phi}$ is any sufficiently regular approximation (in, say, $H^2$), close to $\phi$, that lives in the approximate domain $\hat{\Omega}$ with a sufficiently smooth approximate free boundary $\hat{\Gamma}$ (say, $C^{1,1}$) induced by the approximation $\hat{\boldsymbol\theta}$.
\par
A~key strategy in the derivation of our linearisation consists of the use of higher-order corrections to arrive at more convenient expressions: In particular, since $\hat\phi$ is assumed to be close to~$\phi$, we will often use that $\hat\phi$ satisfies the boundary conditions up to, say, $O(\|\phi-\hat\phi\|, \|\boldsymbol{\theta} - \hat{\boldsymbol{\theta}}\|)$.

\subsection{Linearisation of \texorpdfstring{$\mathcal{R}_1$}{R1}} \label{5.1}
The G$\hat{\text{a}}$teaux derivative at $\hat{\phi}$ in the direction $\delta \phi$ can be evaluated as
\begin{eqnarray}
\left<\partial_{\phi}\mathcal{R}_1\left(\left(\hat{\boldsymbol\theta}, \hat{\phi}\right); v\right), \delta\phi\right> &=& \lim_{t\to0}\frac{\mathcal{R}_1\left(\left(\boldsymbol\theta, \hat{\phi}+t\delta\phi\right);v\right) - \mathcal{R}_1\left(\left(\boldsymbol\theta, \hat{\phi}\right);v\right) }{t}\nonumber\\
&=&  \int_{\hat{\Omega}}\nabla\delta\phi\cdot\nabla v d \Omega +\int_{\partial\Omega\setminus\hat{\Gamma}}\omega\delta\phi vd\Gamma.\label{R1_p}
\end{eqnarray}

Then the linearisation with respect to $\boldsymbol\theta$ can be obtained by applying Hadamard formulas from Theorem \ref{Shape_domain} to \eqref{weak form 1}, assuming $f\in H^1$, which yields
\begin{equation}
\left<\partial_{\boldsymbol{\theta  }}\mathcal{R}_1\left(\left(\hat{\boldsymbol\theta}, \hat{\phi}\right); v\right), \delta\boldsymbol\theta\right> =  \int_{\hat{\Gamma}}\nabla\hat{\phi}\cdot\nabla v \delta\boldsymbol\theta\cdot\boldsymbol n  d \Gamma - \int_{\hat{\Gamma}}fv\delta\boldsymbol\theta\cdot\boldsymbol n d \Gamma.\label{L1}
\end{equation}

The tangential gradient $\nabla_{\Gamma}$ and tangential divergence $\text{div}_{\Gamma}$ satisfy
\begin{equation}
\nabla_{\Gamma}(\cdot) =\nabla(\cdot) - \partial_{\boldsymbol n}(\cdot)\boldsymbol n, \qquad \text{div}_{\Gamma}(\cdot)  = \text{div}(\cdot) - \partial_n(\cdot)\mathbf{n}.\label{tangential_gd}
\end{equation}
By substituting \eqref{tangential_gd} into \eqref{L1} and applying the tangential Green's identity \cite{Shapes, sokolowski1992introduction}, \eqref{L1} can be approximated as \footnote{The integral term over $\partial\Gamma$ is missing in the formula in paper \cite{van2013shape}.}
\begin{eqnarray}
\left<\partial_{\boldsymbol\theta}\mathcal{R}_1\left(\left(\hat{\boldsymbol\theta}, \hat{\phi}\right); v\right), \delta\boldsymbol\theta\right> = \int_{\hat{\Gamma}}\left(\nabla_{\Gamma}\hat{\phi}\cdot\nabla_{\Gamma}v + \partial_{\boldsymbol n}\hat{\phi}\partial_{\boldsymbol n}v\right)\delta\boldsymbol\theta\cdot\mathbf{n} d \Gamma- \int_{\hat{\Gamma}}fv\delta\boldsymbol\theta\cdot\boldsymbol n d \Gamma\nonumber\\
  \approx  - \int_{\hat{\Gamma}}\text{div}_{\Gamma}\left(\delta\boldsymbol\theta\cdot\mathbf{n} \nabla_{\Gamma}\hat{\phi}\right) v  d \Gamma + \int_{\partial\hat{\Gamma}}\left(\delta\boldsymbol\theta\cdot\mathbf{n} \nabla_{\Gamma}\hat{\phi}\cdot\boldsymbol{\tau}\right)vd\Gamma- \int_{\hat{\Gamma}}fv\delta\boldsymbol\theta\cdot\boldsymbol n d \Gamma
\end{eqnarray}
where, due to the Neumann boundary condition \eqref{FB_N} (or \eqref{FBD_N}) and $\hat{\phi}$ being close to $\phi^*$, the related term~$\partial_{\boldsymbol n}\hat{\phi}\delta\boldsymbol \theta\cdot\boldsymbol n$ is of higher order, hence it was neglected. We now use the announced compatibility conditions from Remark \ref{remark_compatibility1}-\ref{remark_compatibility3} to remove the integral over $\partial\hat{\Gamma}$,
\begin{equation}
    \left<\partial_{\boldsymbol\theta}\mathcal{R}_1\left(\left(\hat{\boldsymbol\theta}, \hat{\phi}\right); v\right), \delta\boldsymbol\theta\right>\approx  - \int_{\hat{\Gamma}}\text{div}_{\Gamma}\left(\delta\boldsymbol\theta\cdot\mathbf{n} \nabla_{\Gamma}\hat{\phi}\right) v  d - \int_{\hat{\Gamma}}fv\delta\boldsymbol\theta\cdot\boldsymbol n d \Gamma
    \label{R1_t}
\end{equation}

For Dirichlet condition case, \eqref{R1_t} can be written as
\begin{eqnarray}  \left<\partial_{\boldsymbol\theta}\mathcal{R}_1\left(\left(\hat{\boldsymbol\theta}, \hat{\phi}\right); v\right), \delta\boldsymbol\theta\right> \approx &-& \int_{\hat{\Gamma}}\text{div}_{\Gamma}\left(\delta\boldsymbol\theta\cdot\mathbf{n} \nabla_{\Gamma}h\right) v  d\Gamma- \int_{\hat{\Gamma}}fv\delta\boldsymbol\theta\cdot\boldsymbol n d \Gamma\label{R1_t_D}
\end{eqnarray} 
due to $\hat{\phi} = h + O(||\phi-\hat{\phi}||,||\boldsymbol{\theta}-\hat{\boldsymbol{\theta}}||)$ on the free boundary.

\subsection{Linearisation of \texorpdfstring{$\mathcal{R}_2$}{R2} with Dirichlet condition}

Considering first the Dirichlet boundary condition, we have
\begin{equation}
\mathcal{R}_2\left(\left(\boldsymbol\theta,\phi\right); w\right) = \int_{\Gamma_F}\left(\phi-h\right)w d \Gamma.\label{L2_D}
\end{equation}

Similar to the linearisation of $\mathcal{R}_1$ with respect to $\phi$, it is straightforward to evaluate the G$\hat{\text{a}}$teaux derivative at $\phi$ in the direction $\delta\phi$,
\begin{equation}
\left<\partial_{\boldsymbol\phi}\mathcal{R}_2\left(\left(\hat{\boldsymbol\theta}, \hat{\phi}\right); w\right), \delta\phi\right> = \int_{\hat{\Gamma}}\delta\phi w d \Gamma.\label{R2D_p}
\end{equation}

Then by using the Hadamard formula on the boundary integral \eqref{L2_D}, assuming $h\in H^2$ (and recall that $\hat{\phi}\in H^2$), we have the shape derivative
\begin{eqnarray}
\left<\partial_{\boldsymbol\theta}\mathcal{R}_2\left(\left(\hat{\boldsymbol\theta}, \hat{\phi}\right); w\right), \delta\boldsymbol\theta\right> &=& \int_{\hat{\Gamma}}\left(\partial_{\boldsymbol n}+\kappa\right)\left[\left(\hat{\phi}-h\right) w\right]\delta\boldsymbol\theta\cdot\boldsymbol n  d \Gamma,\nonumber\\
 &=& \int_{\hat{\Gamma}}\left[\partial_{\boldsymbol n}\left(\hat{\phi}-h\right)w + \left(\hat{\phi}-h\right)\partial_{\boldsymbol n}w + \kappa\left(\hat{\phi}-h\right)w\right]\delta\boldsymbol\theta\cdot\boldsymbol n d \Gamma.\nonumber\\ \label{R2_H}
\end{eqnarray}

Using the Dirichlet condition \eqref{FBD_D} and Neumann condition \eqref{FBD_N} on the free boundary, we can neglect the $(\hat{\phi}-h)-$term and $(\partial_{\boldsymbol n}\hat{\phi})-$term in~\eqref{R2_H}, similar to the reasoning in Section~\ref{5.1}. We then have the approximation
\begin{equation}
\left<\partial_{\boldsymbol\theta}\mathcal{R}_2\left(\left(\hat{\boldsymbol\theta}, \hat{\phi}\right); w\right), \delta\boldsymbol\theta\right>  \approx -\int_{\hat{\Gamma}}(\partial_{\boldsymbol n}h) w\delta\boldsymbol\theta\cdot\boldsymbol n d \Gamma.\label{R2D_t}
\end{equation}

\subsection{Linearisation of \texorpdfstring{$\mathcal{R}_2$}{R2} with Bernoulli condition} \label{linearisation_bernoulli}
To perform the linearisation of the Bernoulli condition, we require more regularity on $\nabla\hat\phi$ as well as the test function~$w$. It is sufficient to assume $\hat\phi\in H^3$ and $w\in W^{2,\infty} $.\footnote{The end result \eqref{L2_BC_1} and \eqref{general_deriviation} of the linearisation indicates that these regularity requirements may be weakened, although this has not been pursued further.} 

Substituting the Bernoulli condition \eqref{FB_BC} into the weak form \eqref{weak form 2}, we have
\begin{equation}
\mathcal{R}_2\left(\left(\boldsymbol\theta,\phi\right); w\right) = \int_{\Gamma_F}\left(a\left|\nabla\phi\right|^2 + b x_N + c\right) w d \Gamma.\label{L2_BC}
\end{equation}

The linearisation in terms of $\phi$ at approximation $\hat{\phi}$ is
\begin{equation}
\label{R2B_p}
\left<\partial_{\phi}\mathcal{R}_2\left(\left(\hat{\boldsymbol\theta}, \hat{\phi}\right); w\right), \delta\phi\right> = \int_{\hat{\Gamma}}2a\nabla\hat{\phi}\cdot\nabla\delta\phi w d \Gamma  {}\approx{}\int_{\hat{\Gamma}}2a\nabla_{\Gamma}\hat{\phi}\cdot\nabla_{\Gamma}\delta\phi wd\Gamma, 
\end{equation}
where the normal component was neglected, similar to the reasoning in Section~\ref{5.1}.

To find the G$\hat{\text{a}}$teaux derivative with respect to $\boldsymbol\theta$ at $\hat{\boldsymbol \theta}$, the Hadamard formula yields
\begin{eqnarray}
\left<\partial_{\boldsymbol\theta}\mathcal{R}_2\left(\left(\hat{\boldsymbol\theta}, \hat{\phi}\right); w\right), \delta\boldsymbol\theta\right> = \int_{\hat{\Gamma}}\left(\partial_{\boldsymbol n}+\kappa\right)\left[\left(a\left|\nabla\hat{\phi}\right|^2+b\hat{x}_N+c\right) w\right]\delta\boldsymbol\theta\cdot\boldsymbol n  d \Gamma.\nonumber
\end{eqnarray}

According to the Bernoulli condition \eqref{FB_BC} and $\hat{\phi}$ being close to $\phi$, $a\left|\nabla\hat{\phi}\right|^2+b\hat{x}_N+c$ is close to $0$, similar to what we did in Section \ref{5.1}, the approximation is therefore
\begin{eqnarray}
\left<\partial_{\boldsymbol\theta}\mathcal{R}_2\left(\left(\hat{\boldsymbol\theta}, \hat{\phi}\right); w\right), \delta\boldsymbol\theta\right> \approx \int_{\hat{\Gamma}}\left(a\partial_{\boldsymbol n}\left(\left|\nabla\hat{\phi}\right|^2\right) + b n_N\right)w\delta\boldsymbol\theta\cdot\boldsymbol n d \Gamma, \label{L2_BC_1}
\end{eqnarray}
where $n_N$ is the $x_N$-coordinate (vertical coordinate) of the unit normal vector $\boldsymbol n$. In the two-dimensional case, the value of $n_N$ corresponds to the $y$-component of $\boldsymbol n$. However, in $N$ dimensions, $n_N$ represents the $N$-th component.

We now look more closely at the term $\partial_{\boldsymbol n}\left|\nabla\hat{\phi}\right|^2$.

\subsubsection{\texorpdfstring{$N$}{N} dimensional case}
We first continue assuming the general case in $N$ dimensions, and we will look into the two-dimensional case later for convenience to the reader. 

We introduce the index form of $\nabla$ by
\begin{eqnarray}
    \left(\nabla\hat{\phi}\right)_i = \frac{\partial\hat{\phi}}{\partial x_i} = \hat{\phi}_{,i}, \quad i = 1,\ldots,N,\label{gradient_index}
\end{eqnarray}
such that the Neumann boundary condition \eqref{FB_N} can be rewritten in the form 
\begin{equation}
    \partial_{\boldsymbol n}\hat{\phi} = n_i \hat{\phi}_{,i} = 0,\label{Neumann_index_form}
\end{equation}
where we employ the Einstein summation convention.

Taking the tangential gradient gives:
\begin{equation}
    \nabla_{\Gamma}\left(\partial_{\boldsymbol n}\hat{\phi}\right) = 0.\label{tangential_gradient_NB}
\end{equation}
We define the tangential gradient and the matrix $\left[\nabla_{\Gamma}\boldsymbol x\right]$ as in \cite{dziuk2013finite} 
\begin{eqnarray}
    \left(\nabla_{\Gamma}\hat{\phi}\right)_{\alpha} = \bar{D}_{\alpha}\hat{\phi} &=& \hat{\phi}_{,\alpha}, \quad \alpha = 1,\ldots, N, \label{tangential_gradient_index}\\
    \left[\nabla_{\Gamma}\boldsymbol x\right]_{\alpha i} &=& \bar{D}_{\alpha}x_i = x_{i,\alpha}.\label{tangential_gradient_matrix}
\end{eqnarray}

According to the definition of tangential gradient, we have
\begin{eqnarray}
    \left[\nabla_{\Gamma}\nabla\hat{\phi}\right] &=& \left[\nabla\nabla\hat{\phi}\right] - \partial_{\boldsymbol n}\left(\nabla\hat{\phi}\right)\boldsymbol n^T\nonumber\\
    &=& \left[\nabla\nabla_{\Gamma}\hat{\phi}\right] - \partial_{\boldsymbol n}\left(\nabla_{\Gamma}\hat{\phi}\right)\boldsymbol n^T\label{derivative_assumption}
\end{eqnarray}
where the second step is obtained by using Neumann boundary condition \eqref{FB_N}. Since the $\nabla_{\Gamma}\hat{\phi}$ is the tangential gradient of $\hat{\phi}$ which is only defined on the free surface $\hat{\Gamma}$, it can be extended as a constant beyond the surface such that its normal derivative is zero. Hence,
\begin{equation}
    \left[\nabla_{\Gamma}\nabla\hat{\phi}\right] = \left[\nabla\nabla_{\Gamma}\hat{\phi}\right].\label{tangential_gradient_gradient_assumption}
\end{equation}

The equation \eqref{tangential_gradient_NB} is equivalent to
\begin{eqnarray}
     \left(\nabla_{\Gamma}\left(\partial_{\boldsymbol n}\hat{\phi}\right)\right)_{\alpha}
 &=&\bar{D}_{\alpha}\left(n_i\hat{\phi}_{,i}\right) \nonumber\\
 &=& n_{i,\alpha}\hat{\phi}_{,i} + n_i\hat{\phi}_{,i\alpha} \nonumber\\
 &=& \left(\left[\nabla_{\Gamma}\boldsymbol n\right]\,\nabla\hat{\phi} +\left[\nabla_{\Gamma}\nabla\hat{\phi}\right]\, \boldsymbol n\right)_{\alpha}\nonumber\\
 &=& 0.\label{tangential_gradient_NB_index}
\end{eqnarray}

By using \eqref{tangential_gradient_gradient_assumption} and \eqref{tangential_gradient_NB_index}, we have
\begin{alignat}{2}
    \partial_{\boldsymbol n}\left(\left|\nabla\hat{\phi}\right|^2\right) &= \partial_{\boldsymbol n}\left(\left|\nabla_{\Gamma}\hat{\phi}\right|^2\right)\nonumber\\
    &=\partial_{\boldsymbol n}\left(\hat{\phi}_{,\alpha}\,\hat{\phi}_{,\alpha}\right)\nonumber\\
    &= 2n_i\hat{\phi}_{,i\alpha}\hat{\phi}_{,\alpha}\tag{by \eqref{tangential_gradient_gradient_assumption}}\\
    &= 2\left(\nabla_{\Gamma}\hat{\phi}\right)^T\cdot\left[\nabla_{\Gamma}\nabla\hat{\phi}\right]\cdot\boldsymbol n\nonumber\\
    &= -2\left(\nabla_{\Gamma}\hat{\phi}\right)^T\cdot\left[\nabla_{\Gamma}\boldsymbol n\right]\cdot\nabla_{\Gamma}\hat{\phi},
    \label{general_deriviation}
\end{alignat}
where $\left[\nabla_{\Gamma}\boldsymbol n\right]$ is the extended Weingarten map \cite{dziuk2013finite}, which is a tensor containing curvature type quantities. In particular, the trace of $\left[\nabla_{\Gamma}\boldsymbol n\right]$ coincides with the summed curvature \cite[Sec. 4.5.2]{walker2015shapes}.

Substituting~\eqref{general_deriviation} into \eqref{L2_BC_1}, the (approximate) shape linearisation in the $N$-dimensional case        becomes
\begin{equation}
\left<\partial_{\boldsymbol\theta}\mathcal{R}_2\left(\left(\hat{\boldsymbol\theta}, \hat{\phi}\right); w\right), \delta\boldsymbol\theta\right> \approx \int_{\hat{\Gamma}}\left[-2a\left(\nabla_{\Gamma}\hat{\phi}\right)^T\cdot\left[\nabla_{\Gamma}\boldsymbol n\right]\cdot\nabla_{\Gamma}\hat{\phi} + b n_N\right]w\delta\boldsymbol\theta\cdot\boldsymbol n d \Gamma.\label{3d_linearisation}
\end{equation}

\subsubsection{Three dimensional case}\label{section_3d}
For the three-dimensional case, let $\kappa_1$ and $\kappa_2$ be the principle curvatures. The matrix~$\left[\nabla_{\Gamma}\boldsymbol n\right]$ has eigenvalues $\{0, \kappa_1, \kappa_2\}$ and corresponding normalised eigenvectors $\{\boldsymbol{0}, \boldsymbol{d}_1, \boldsymbol{d}_2\}$ \cite[Sec. 4.5.2]{walker2015shapes}. Since $\left[\nabla_{\Gamma}\boldsymbol n\right]$ is symmetric \cite{walker2015shapes, dziuk2013finite}, by the spectral decomposition theorem, 
$$\left[\nabla_{\Gamma}\boldsymbol n\right] = Q\Lambda Q^{T}$$
where $$\Lambda = \begin{pmatrix}
    0 & 0 & 0 \\
    0 & \kappa_1 & 0  \\
    0 & 0 & \kappa_2 
\end{pmatrix}, \quad Q = \begin{pmatrix}
    \boldsymbol{0} & \boldsymbol{d}_1 & \boldsymbol{d}_2
\end{pmatrix}.$$
Hence,
\begin{eqnarray}
\left(\nabla_{\Gamma}\hat{\phi}\right)^T\cdot\left[\nabla_{\Gamma}\boldsymbol n\right]\cdot\nabla_{\Gamma}\hat{\phi} &=& \left(\nabla_{\Gamma}\hat{\phi}\right)^T\cdot Q \Lambda Q^{T}\cdot\nabla_{\Gamma}\hat{\phi}\nonumber\\
&=& \begin{pmatrix}
    0 & \boldsymbol{d}_1\cdot\nabla_{\Gamma}\hat{\phi} & \boldsymbol{d}_2\cdot\nabla_{\Gamma}\hat{\phi}
\end{pmatrix}\cdot\Lambda\cdot\begin{pmatrix}
    0 \\ \boldsymbol{d}_1\cdot\nabla_{\Gamma}\hat{\phi} \\ \boldsymbol{d}_2\cdot\nabla_{\Gamma}\hat{\phi}
\end{pmatrix}\nonumber\\
&=& \kappa_1(\boldsymbol{d}_1\cdot\nabla_{\Gamma}\hat{\phi})^2 + \kappa_2(\boldsymbol{d}_2\cdot\nabla_{\Gamma}\hat{\phi})^2\label{3d}
\end{eqnarray}
where the last step is obtained because $\boldsymbol{d}_1$ and $\boldsymbol{d}_2$ are orthonormal \cite[Sec. 3.2.4]{walker2015shapes}.

Substituting~\eqref{3d} into~\eqref{L2_BC_1}, the approximate shape linearisation (up to higher-order terms) in the 3-D case becomes
\begin{equation}
\left<\partial_{\boldsymbol\theta}\mathcal{R}_2\left(\left(\hat{\boldsymbol\theta}, \hat{\phi}\right); w\right), \delta\boldsymbol\theta\right> \approx \int_{\hat{\Gamma}}[-2a\kappa_1(\boldsymbol{d}_1\cdot\nabla_{\Gamma}\hat{\phi})^2 -2a \kappa_2(\boldsymbol{d}_2\cdot\nabla_{\Gamma}\hat{\phi})^2 + b n_z]w\delta\boldsymbol\theta\cdot\boldsymbol n d \Gamma,\label{Nd_linearisation}
\end{equation}
where $n_z$ is the $z$-component of the unit normal vector $\boldsymbol{n}$.

\subsubsection{Two dimensional case}
For convenience to the reader, now we look into a specific case of the linerisation of $\mathcal{R}_2$, namely in two dimensions using Cartesian coordinates, which is more direct. In this case, we assume that we can introduce $\eta(x)$ as the vertical displacement of the free surface with respect to the referenced free surface, the horizontal $x$-axis, such that $\boldsymbol \theta = \left(x, \eta(x)\right)$. 

Given the approximation $\hat{\boldsymbol\theta} = \left(x, \hat{ \eta}(x)\right)$, we have the unit normal vector $\boldsymbol n = \frac{1}{\sqrt{1+\hat{ \eta}_x^2}}\left(-\hat{ \eta}_x,1\right)$ and the unit tangential vector $\boldsymbol\tau = \frac{1}{\sqrt{1+\hat{ \eta}_x^2}}\left(1, \hat{ \eta}_x\right)$. Then the Neumann boundary condition \eqref{FB_N} on the free boundary can be written in the form of 
\begin{equation}
-\hat{ \eta}_x\hat{\phi}_x + \hat{\phi}_y = 0.\nonumber
\end{equation}
This implies that its tangential derivative is also zero, i.e.
\begin{equation}
\left(\boldsymbol\tau\cdot\nabla\right)\left(-\hat{ \eta}_x\hat{\phi}_x + \hat{\phi}_y\right) = 0,\nonumber
\end{equation}
which is equivalent to
\begin{equation}
-\hat{ \eta}_{xx}\hat{\phi}_x - \hat{ \eta}_x\hat{\phi}_{xx}+\hat{\phi}_{xy}-\hat{ \eta}_x^2\hat{\phi}_{xy}+\hat{ \eta}_x\hat{\phi}_{yy}=0.\label{TD}
\end{equation}

Then we have
\begin{eqnarray}
\partial_{\boldsymbol n}\left(\left|\nabla\hat{\phi}\right|^2\right) &=& \frac{1}{\sqrt{1+\hat{ \eta}_x^2}}\left(-\hat{ \eta}_x\partial_x + \partial_y\right)\left(\hat{\phi}_x^2+\hat{\phi}_y^2\right)\nonumber\\
 &=& \frac{2}{\sqrt{1+\hat{ \eta}_x^2}}\hat{\phi}_x\left(-\hat{ \eta}_x\hat{\phi}_{xx}-\hat{ \eta}_x^2\hat{\phi}_{xy}+\hat{\phi}_{xy}+\hat{ \eta}_x\hat{\phi}_{yy}\right)\nonumber\\
 &=&\frac{2}{\sqrt{1+\hat{ \eta}_x^2}}\hat{ \eta}_{xx}\left(\hat{\phi}_x\right)^2\nonumber\\
 &=&2\kappa\left|\nabla\hat{\phi}\right|^2, \label{dn}
\end{eqnarray}
where $\kappa = \partial_x\left(\frac{\hat{ \eta}_x}{\sqrt{1+\hat{ \eta}_x^2}}\right)$, which is the curvature. The second and last steps are obtained by substituting the Neumann condition, and the third step is obtained by substitution of \eqref{TD}.

In the two dimensional case, we have
\begin{eqnarray}
    \left[\nabla_{\Gamma}\boldsymbol n\right] &=& \left[\nabla\boldsymbol n\right] - \partial_{\boldsymbol n}\left(\boldsymbol n\right)\boldsymbol n^T\nonumber\\
    &=& \begin{pmatrix}
        -\kappa+\kappa n_x & -\eta_x\kappa+\kappa n_x n_y\\
        \kappa\eta_x n_x^2 & \kappa\eta_x n_x n_y
    \end{pmatrix}\label{matrix_2d}
\end{eqnarray}
where $\boldsymbol n = \left(n_x, n_y\right) = \frac{1}{\sqrt{1+\hat{ \eta}_x^2}}\left(-\hat{ \eta}_x,1\right).$ By substituting \eqref{matrix_2d} into \eqref{general_deriviation} and using the Neumann condition \eqref{FB_N}, \eqref{dn} is consistent with \eqref{general_deriviation}. The details can be found in Appendix \ref{Appendix: calculation}.

On substitution from \eqref{dn} into \eqref{L2_BC_1}, the approximate shape linearisation in the 2-D case        is
\begin{equation}
\left<\partial_{\boldsymbol\theta}\mathcal{R}_2\left(\left(\hat{\boldsymbol\theta}, \hat{\phi}\right); v\right), \delta\boldsymbol\theta\right> \approx \int_{\hat{\Gamma}}\left(2a\kappa\left|\nabla\hat{\phi}\right|^2 + b n_y\right)w\delta\boldsymbol\theta\cdot\boldsymbol n d \Gamma.\label{R2B_t}
\end{equation}

\section{Newton-Like Schemes}\label{chap:3.6} 

Next, we use the linearisations in the previous section to construct Newton-like schemes. We first consider the general case in $\mathbb{R}^N$. \par
We introduce $\boldsymbol\theta = \hat{\boldsymbol\theta} + \delta\boldsymbol\theta \in C^{0,1}(\Gamma_0,\mathbb{R}^N)$ and $\phi = \hat{\phi} + \delta\phi \in H^1(\mathbb{R}^N)$, where  $\delta\boldsymbol\theta\in C^{0,1}(\Gamma_0;\mathbb{R}^N)$ and $\delta\phi \in H^1(\hat{\Omega})\subset H^1(\mathbb{R}^N)$ are the corrections to $\hat{\boldsymbol{\theta}} \in C^{0,1}(\Gamma_0;\mathbb{R}^N)$, which generates the domain~$\hat{\Omega}$ with free boundary~$\hat{\Gamma} = \Gamma_0 + \hat{\boldsymbol{\theta}}$, and $\hat{\phi}\in H^1(\mathbb{R}^N)$, respectively.%
\footnote{The inclusion $H^1(\hat{\Omega})\subset H^1(\mathbb{R}^N)$ is meant in the sense that each $\delta\phi \in H^1(\hat{\Omega})$ has a (non-unique) extension onto~$\mathbb{R}^N \setminus \hat{\Omega}$, which is in $H^1(\mathbb{R}^N)$.}
\par
In each iteration, a reference free boundary~$\hat{\Gamma}$ is updated, and thereby the reference domain~$\hat{\Omega}$. The exact Newton method for $\left(\delta\boldsymbol\theta, \delta\phi\right)$, in weak form, would be 
\begin{subequations}
\begin{eqnarray}
\left<\partial_{(\boldsymbol\theta, \phi)}\mathcal{R}_1\left(\left(\hat{\boldsymbol\theta}, \hat{\phi}\right);v\right), \left(\delta\boldsymbol\theta, \delta\phi\right)\right> &=& - \mathcal{R}_1\left(\left(\hat{\boldsymbol\theta}, \hat{\phi}\right); v\right) \quad \forall v\in V,\label{N_1} \\
\left<\partial_{(\boldsymbol\theta, \phi)}\mathcal{R}_2\left(\left(\hat{\boldsymbol\theta}, \hat{\phi}\right);w\right), \left(\delta\boldsymbol\theta, \delta\phi\right)\right> &=& - \mathcal{R}_2\left(\left(\hat{\boldsymbol\theta}, \hat{\phi}\right); v\right) \quad \forall w\in W.\label{N_2}
\end{eqnarray}
\end{subequations}
Instead, we obtain more convenient Newton-like schemes by using the higher-order corrections of Section~\ref{chap:3.5} to the above derivatives.%
\footnote{In particular, when $\hat{\phi}$ and $\hat{\boldsymbol\theta}$ are the exact solutions $(\boldsymbol{\theta}^*, \phi^*)$}, the Newton-like schemes coincides with the exact Newton scheme.
We subsequently provide a strong form interpretation of the scheme.

\subsection{Weak form of the problem with Dirichlet Boundary condition}
The Newton-like equation for $\mathcal{R}_1$ is obtained by combining \eqref{R1_p} and the approximation \eqref{R1_t_D} of $\partial_{\boldsymbol\theta}\mathcal{R}_1\left(\left(\hat{\boldsymbol\theta}, \hat{\phi}\right); v\right)$, i.e.,
\begin{subequations}
\begin{eqnarray}
\int_{\hat{\Omega}}\nabla\delta\phi\cdot\nabla v d \Omega + \int_{\partial\hat{\Omega}}\setminus\hat{\Gamma}\omega\delta\phi vd\Gamma- \int_{\hat{\Gamma}}\text{div}_{\Gamma}\left(\delta\boldsymbol\theta\cdot\boldsymbol n\nabla_{\Gamma}h\right)v d \Gamma - \int_{\hat{\Gamma}}fv\delta\boldsymbol\theta\cdot\boldsymbol n d \Gamma
\nonumber
\\ = -\mathcal{R}_1\left(\left(\hat{\boldsymbol\theta}, \hat{\phi}\right); v\right), \quad \forall v\in V.
\label{R1_NL_D}
\end{eqnarray}

For the Dirichlet boundary condition, the Newton-like equation for~$\mathcal{R}_2$ is derived based on \eqref{R2D_p} and approximation \eqref{R2D_t} as
\begin{eqnarray}
\int_{\hat{\Gamma}}\delta\phi w d \Gamma - \int_{\hat{\Gamma}}\left(\partial_{\boldsymbol n}h\right)w\delta\boldsymbol\theta\cdot\boldsymbol n d \Gamma = -\mathcal{R}_2\left(\left(\hat{\boldsymbol\theta}, \hat{\phi}\right); w\right), \quad \forall w\in W.
\label{R2D_NL}
\end{eqnarray}
\end{subequations}

\subsection{Weak form of the problem with Bernoulli Boundary condition}
The Newton-like equation for $\mathcal{R}_1$ is obtained by combining \eqref{R1_p} and the approximation \eqref{R1_t} of $\partial_{\boldsymbol\theta}\mathcal{R}_1\left(\left(\hat{\boldsymbol\theta}, \hat{\phi}\right); v\right)$, i.e.,
\begin{subequations}
\begin{eqnarray}
\int_{\hat{\Omega}}\nabla\delta\phi\cdot\nabla v d \Omega + \int_{\partial\hat{\Omega}\setminus\hat{\Gamma}}\omega\delta\phi v d\Gamma- \int_{\hat{\Gamma}}\text{div}_{\Gamma}\left(\delta\boldsymbol\theta\cdot\boldsymbol n\nabla_{\Gamma}\hat{\phi}\right)v d \Gamma - \int_{\hat{\Gamma}}fv\delta\boldsymbol\theta\cdot\boldsymbol n d \Gamma
\nonumber
\\ = -\mathcal{R}_1\left(\left(\hat{\boldsymbol\theta}, \hat{\phi}\right); v\right), \quad \forall v\in V.
\label{R1_NL}
\end{eqnarray}

For the Bernoulli condition, introducing \eqref{R2B_p}, \eqref{L2_BC_1} and \eqref{general_deriviation}, the Newton-like equation for~$\mathcal{R}_2$ is 
\begin{eqnarray}
\int_{\hat{\Gamma}}2a\nabla_{\Gamma}\hat{\phi}\cdot\nabla_{\Gamma}\delta\phi\, w\, d \Gamma + \int_{\hat{\Gamma}}\left(-2\left(\nabla_{\Gamma}\hat{\phi}\right)^T\cdot\left[\nabla_{\Gamma}\boldsymbol n\right]\cdot\nabla_{\Gamma}\hat{\phi}+bn_N\right)w\, \delta\boldsymbol\theta\cdot\boldsymbol n \,d \Gamma 
\nonumber
\\ = -\mathcal{R}_2\left(\left(\hat{\boldsymbol\theta}, \hat{\phi}\right); w\right), \quad \forall w\in W.
\label{R2B_NL}
\end{eqnarray}
\end{subequations}

\subsection{Strong form: General free-boundary perturbations}

\begin{table}[t!]
\footnotesize
\caption{The coupled shape-Newton scheme solving for $\left(\delta\boldsymbol\theta, \delta\phi\right)$ using a linearised Bernoulli boundary condition~\eqref{table_1_Bernoulli} on the free boundary. For the linearised Dirichlet boundary condition on the free boundary, replace~\eqref{table_1_Bernoulli} by~\eqref{CS_D}.}
	\fbox{%
		\parbox{0.95\linewidth}{%
			\begin{enumerate}
				\item Initialize with $\left(\boldsymbol\theta^0, \phi^0\right)$; set $k = 0$.\\
				\item  Given $\big(\boldsymbol\theta^k, \phi^k\big)$, solve the linear coupled 
    problem for $\left(\delta\boldsymbol\theta\cdot\boldsymbol n, \delta\phi\right)$:  
    \begin{subequations}
    \begin{eqnarray}
        \nabla^2\delta\phi = -\nabla^2\phi^k-f\quad \text{in }\Omega^k,\label{table1_laplace}\\
\partial_n\delta\phi  \,-\,\text{div}_{\Gamma}\left(\delta\boldsymbol\theta\cdot\mathbf{n}\nabla_{\Gamma}\phi^k\right) - f\delta\boldsymbol\theta\cdot\boldsymbol n = -\partial_{\boldsymbol{n}}\phi^k, \quad \text{on }\Gamma^k,
\\
\partial_{\boldsymbol n}\delta\phi + \omega\delta\phi = g+\omega h-\left(\partial_{\boldsymbol n}\phi^k + \omega\phi^k\right) \quad \text{on }\partial\Omega^k\setminus\Gamma^k,\label{table1_fixed}
\\
\label{table_1_Bernoulli}
2a\nabla_{\Gamma}\phi^k\cdot\nabla_{\Gamma}\delta\phi + \left(-2\left(\nabla_{\Gamma}\phi^k\right)^T\cdot\left[\nabla_{\Gamma}\boldsymbol n\right]\cdot\nabla_{\Gamma}\phi^k+bn_N\right)\delta\boldsymbol\theta\cdot\boldsymbol n \\
=  -\left(a\left|\nabla\phi^k\right|^2+b\hat{x}_N+c\right), \quad \text{on } \Gamma^k,\nonumber
    \end{eqnarray}
    \end{subequations}
				\item Update the free boundary displacement and potential as
\begin{alignat*}{2}
				\boldsymbol\theta^{k+1} &= \boldsymbol\theta^{k} + \delta\boldsymbol\theta\,,
\\
				\phi^{k+1} &= \phi^k+\delta\phi\,.
				\end{alignat*}
    \item Update the free boundary (hence the domain) as
   \begin{alignat*}{2}
   \Gamma^{k+1} &= \Gamma_0 + \boldsymbol{\theta}^{k+1}
   \end{alignat*}
    		\item[~]
	Then repeat from step~2~with $k:=k+1$ until convergence.
 \end{enumerate}
		}%
		}
	\label{dt}
\end{table}

It is important to provide a strong form interpretation of the Newton-like scheme, so that the linearised equations can be used by methods that don't use weak forms. Furthermore, the strong form provides further insight and a starting point for analysis.

In the Dirichlet case, the strong form problem for $\left(\delta\boldsymbol\theta \cdot \boldsymbol{n}\,,\,\delta\phi\right)$ extracted from \eqref{R1_NL_D}-\eqref{R2D_NL} is:%
\footnote{We note that~\cite[Section~4.1]{van2013shape} has several typos in the strong from of the linearised system. Equations~\eqref{SF_laplace_D}--\eqref{CS_D} are  correct versions for the case in~\cite[Section~4.1]{van2013shape} with vanishing Neumann data (i.e., set their $g=0$).}

\begin{subequations}
\label{linSys:Dir}
\begin{alignat}{2}
\nabla^2\delta\phi &= -\nabla^2\hat{\phi}-f &\quad &\text{in }\hat{\Omega},\label{SF_laplace_D}
\\
\partial_{\boldsymbol n}\delta\phi  + \omega\delta\phi &= g+\omega h-\left(\partial_{\boldsymbol n}\hat{\phi} + \omega\hat{\phi}\right) &\quad & \text{on }\partial\hat{\Omega}\setminus\hat{\Gamma},
\\
\partial_n\delta\phi - \text{div}_{\Gamma}\left(\delta\boldsymbol\theta\cdot\boldsymbol{n}\nabla_{\Gamma}h\right) - f\delta\boldsymbol\theta\cdot\boldsymbol n &=-\partial_{\boldsymbol{n}}\hat{\phi} , &\quad &\text{on }\hat{\Gamma}, 
\label{CS_ND}\\
\delta\phi - \partial_{\boldsymbol n}h\,\delta\boldsymbol\theta\cdot\boldsymbol n &= h-\hat{\phi}, &\quad &\text{on }\hat{\Gamma}
\label{CS_D}
\end{alignat}
\end{subequations}
while in the case of the Bernoulli condition, the strong form problem for $\left(\delta\boldsymbol\theta \cdot \boldsymbol{n}\,,\,\delta\phi\right)$ extracted from~\eqref{R1_NL}-\eqref{R2B_NL} is:
\begin{subequations}
\label{linSys:Bern}
\begin{alignat}{2}
\nabla^2\delta\phi &= -\nabla^2\hat{\phi}-f &\quad &\text{in }\hat{\Omega},
\label{CS_L}
\\
\partial_{\boldsymbol n}\delta\phi  + \omega\delta\phi &= g+\omega h-\left(\partial_{\boldsymbol n}\hat{\phi} + \omega\hat{\phi}\right) &\quad & \text{on }\partial\hat{\Omega}\setminus\hat{\Gamma},
\label{CS_R}
\\
\partial_n\delta\phi {}-{} \text{div}_{\Gamma}\left(\delta\boldsymbol\theta\cdot\mathbf{n}\nabla_{\Gamma}\hat{\phi}\right) - f\delta\boldsymbol\theta\cdot\boldsymbol n &={-}\partial_{\boldsymbol{n}}\hat{\phi} , &\quad  &\text{on }\hat{\Gamma}, 
\label{CS_N}
\\ \notag
2a\nabla_{\Gamma}\hat{\phi}\cdot\nabla_{\Gamma}\delta\phi + \left(-2\left(\nabla_{\Gamma}\hat{\phi}\right)^T\cdot\left[\nabla_{\Gamma}\boldsymbol n\right]\cdot\nabla_{\Gamma}\hat{\phi}+bn_N\right)\delta\boldsymbol\theta\cdot\boldsymbol n \hspace*{-12em}
\\ &= -\left(a\left|\nabla\hat{\phi}\right|^2+b\hat{x}_N+c\right), &\quad  &\text{on } \hat{\Gamma}.\label{CS_B}
\end{alignat}
\end{subequations}

The iterative algorithm associated to~\eqref{CS_L}--\eqref{CS_B} is given in Table~\ref{dt}. The solutions are updated as
 $$\boldsymbol\theta^{k+1} = \boldsymbol\theta^k+\delta\boldsymbol\theta, \quad \phi^{k+1} = \phi^k +\delta\phi,$$ 
where
$\delta\boldsymbol{\theta}$ is such that $\delta\boldsymbol{\theta}\cdot\boldsymbol{n}$ satisfies the problem \eqref{table1_laplace}--\eqref{table_1_Bernoulli} (while its tangential component is free to specify).
Accordingly, the free boundary is updated as $\Gamma^{k+1} = \Gamma_0 + \boldsymbol{\theta}^{k+1}$.

\begin{remark}[Solving directly for~$\phi^{k+1}$] 
One can write the linearized system in mixed total/update form, which solves for the variables~$(\phi^{k+1},\delta\boldsymbol{\theta})$, instead of~$(\delta\phi,\delta\boldsymbol{\theta})$. This can be particularly helpful to remove any dependencies on $\phi^k$ (which lives on the previous domain~$\Omega^{k-1}$, hence would need a suitable extension onto~$\Omega^k$); see~\cite[Remark~5]{van2013shape} where, in case of the Dirichlet boundary condition, the dependence on~$\phi^k$ is shown to be completely eliminated.
\end{remark}

\begin{remark}[Solvability of the shape-linearized systems]
\label{remark_solve_dirichlet}
Both linearized systems~\eqref{linSys:Dir} and~\eqref{linSys:Bern} can be shown to have a unique weak solution under certain conditions of the data. We have presented the details of these well-posedness analyses in Appendix~\ref{app:solvability_dirichlet} and~\ref{app:solvability_bernoulli}, for \eqref{linSys:Dir} and~\eqref{linSys:Bern} respectively. In both cases, the analysis establishes coercivity of a bilinear form for a weak formulation for the variable~$\delta \phi$, obtained by eliminating the variable~$\delta\boldsymbol{\theta}\cdot \boldsymbol{n}$ from the system. 
\par
In the case of system~\eqref{linSys:Dir}, the bilinear form corresponds to that of a Laplacian with a generalized Robin boundary condition involving an oblique derivative. Such problems have been analyzed in, e.g., \cite{rozanov2003analysis, troianiello2013elliptic, wachsmuth2013optimal}. 
\par
In the case of system~\eqref{linSys:Bern}, the bilinear corresponds to that of a Laplacian with a generalized Robin boundary condition involving a surface Laplacian (Laplace--Beltrami operator). Such problems have been analyzed in, e.g., \cite{elliott2013finite, AltAML2019, AltVerMMS2021}. 
\end{remark}

\subsection{Strong form: Vertical free-boundary perturbations}
\begin{table}[t!]
\footnotesize
\caption{The coupled shape-Newton scheme for $(\delta \eta,\delta\phi)$. }
	\fbox{%
		\parbox{0.95\linewidth}{%
			\begin{enumerate}
				\item Initialize with $\left(\eta^0, \phi^0\right)$; set $k = 0$.\\
				\item  Given $\left(\eta^k, \phi^k\right)$, solve the free boundary problem 
    \begin{subequations}
        \begin{gather}
        \nabla^2\delta\phi = -\nabla^2\phi^k-f\quad \text{in }\Omega^k,\label{table2_laplacian}\\
\partial_n\delta\phi+\partial_{\boldsymbol{n}}\hat{\phi} - \text{div}_{\Gamma}\left(\delta\boldsymbol\theta\cdot\mathbf{n}\nabla_{\Gamma}\phi^k\right) - f\delta\boldsymbol\theta\cdot\boldsymbol n = 0, \quad \text{on }\Gamma^k,\label{table2_neumann}\\
\\
2a\nabla_{\Gamma}\phi^k\cdot\nabla_{\Gamma}\delta\phi + \left(2a\kappa\left|\nabla\phi^k\right|^2 + b n_y\right)\delta\eta = -\left(a\left|\nabla\phi^k\right|^2+b\eta^k+c\right), \quad \text{on } \Gamma^k,\nonumber\\
\partial_{\boldsymbol n}\delta\phi + \omega\delta\phi = g+\omega h-\left(\partial_{\boldsymbol n}\phi^k + \omega\phi^k\right), \quad \text{on }\partial\Omega^k\setminus\Gamma^k,\label{table2_fixed}
    \end{gather}
        \end{subequations}
    for $\left(\delta \eta, \delta\phi\right)$, where $\delta \eta = \sqrt{1+ \left(\eta^k_x\right)^2}\left(\delta\boldsymbol\theta\cdot\boldsymbol n\right) = \delta\boldsymbol\theta\cdot\left(-\eta^k_x, 1\right)$.\\
				\item Update the free boundary displacement and potential as
				\begin{alignat*}{2}
				 \eta^{k+1} &=  \eta^{k} + \delta \eta\,,
     \\
\phi^{k+1} &= \phi^k+\delta\phi\,.
				\end{alignat*}
    \item Update the free boundary (hence the domain) as
   \begin{alignat*}{2}
   \Gamma^{k+1} &= \Gamma_0 + (0,\eta^{k+1})
   \end{alignat*}
\item[~]    
    Then repeat from step~2 with $k:=k+1$ until convergence.
			\end{enumerate}
		}%
		}
	\label{de}
\end{table}

A particular scenario arises in a two-dimensional case, where the free boundary is adjusted vertically such that $\boldsymbol{\theta} = (0,\eta)$.
\par
In that case, we have  $ d \Gamma =  d s = \sqrt{1+\hat{ \eta}_x^2} d x$ such that
\begin{equation}
\int_{\hat{\Gamma}}\left(\cdot\right)\delta\boldsymbol\theta\cdot\boldsymbol n d \Gamma = \int_{\hat{\Gamma}}\left(\cdot\right)\delta \eta d x, \label{s_to_x}
\end{equation}
where $s$ is the arc length and $\delta \eta = \sqrt{1+\hat{ \eta}_x^2}\left(\delta\boldsymbol\theta\cdot\boldsymbol n\right) = \delta\boldsymbol\theta\cdot\left(-\hat{ \eta}_x, 1\right)$. The boundary integrals can be evaluated in a referenced domain along the $x$ direction, and this problem can be solved in terms of the pair $(\delta \eta, \delta\phi)$. The algorithm is now displayed as Table \ref{de}, and the geometry is updated vertically with $\delta\eta$. 

\section{Numerical experiments}\label{chap:3.8}
\begin{figure}[t]
  \centering
  \subfloat[The initial domain and the triangulation.]{\includegraphics[width=0.5\linewidth]{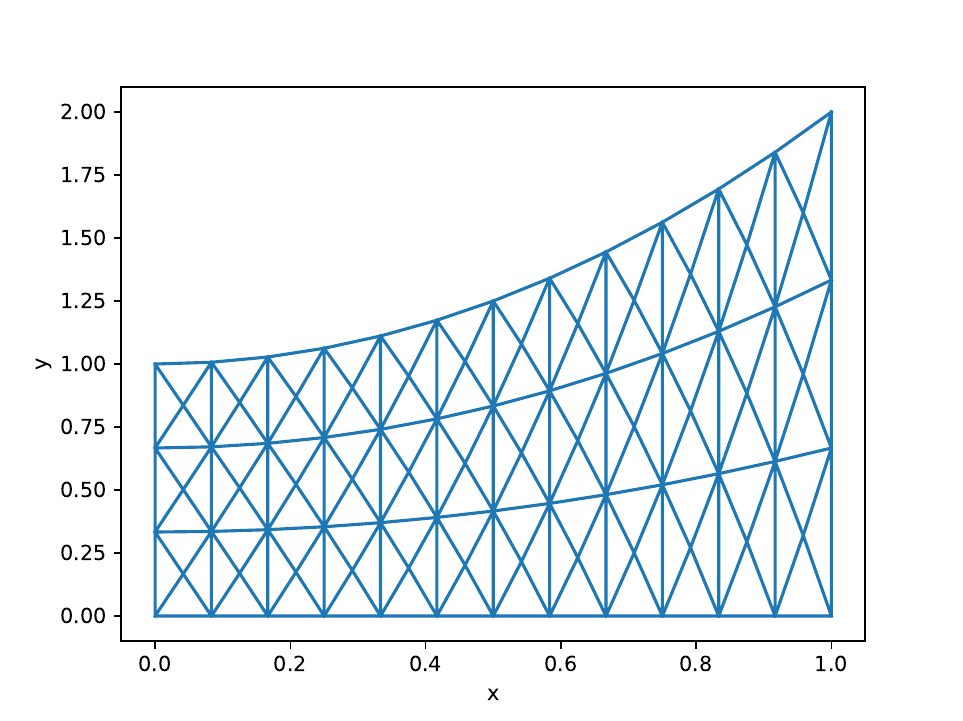}}\hfill
  \subfloat[The domain and the triangulation after the first iteration.]{\includegraphics[width=0.5\linewidth]{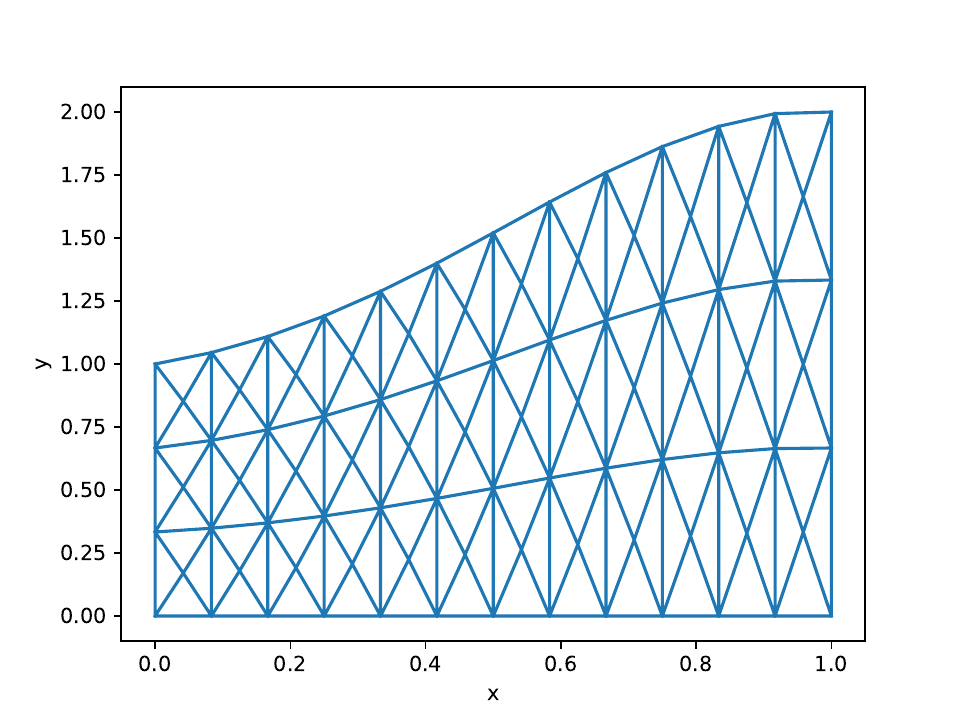}}
  
  \subfloat[The domain and the triangulation after the second iteration.]{\includegraphics[width=0.5\linewidth]{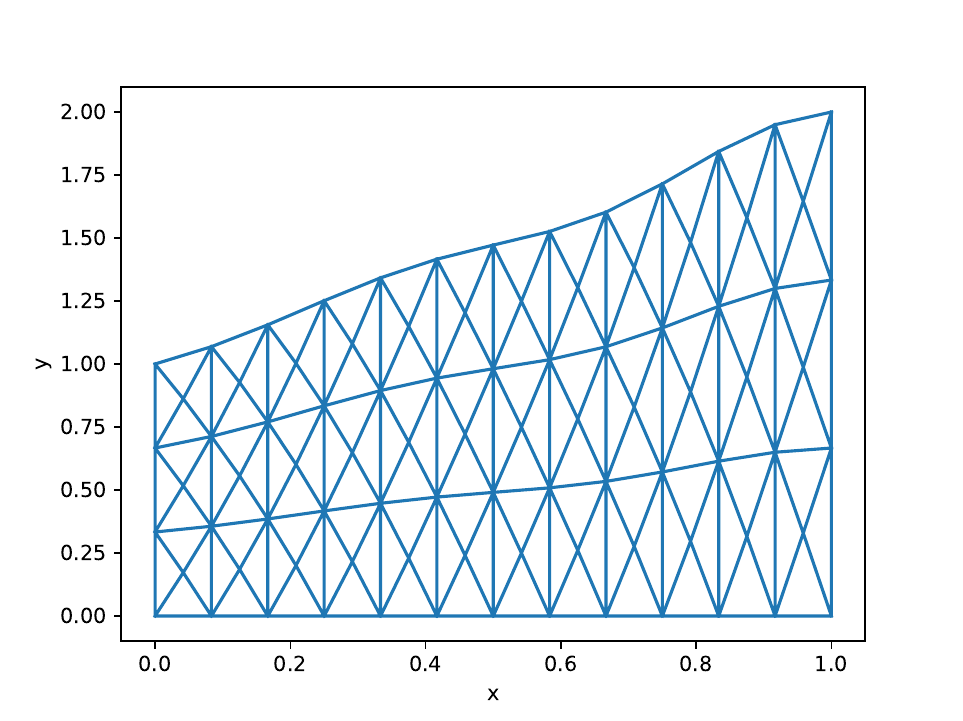}}\hfill
  \subfloat[The domain and the triangulation after the third iteration.]{\includegraphics[width=0.5\linewidth]{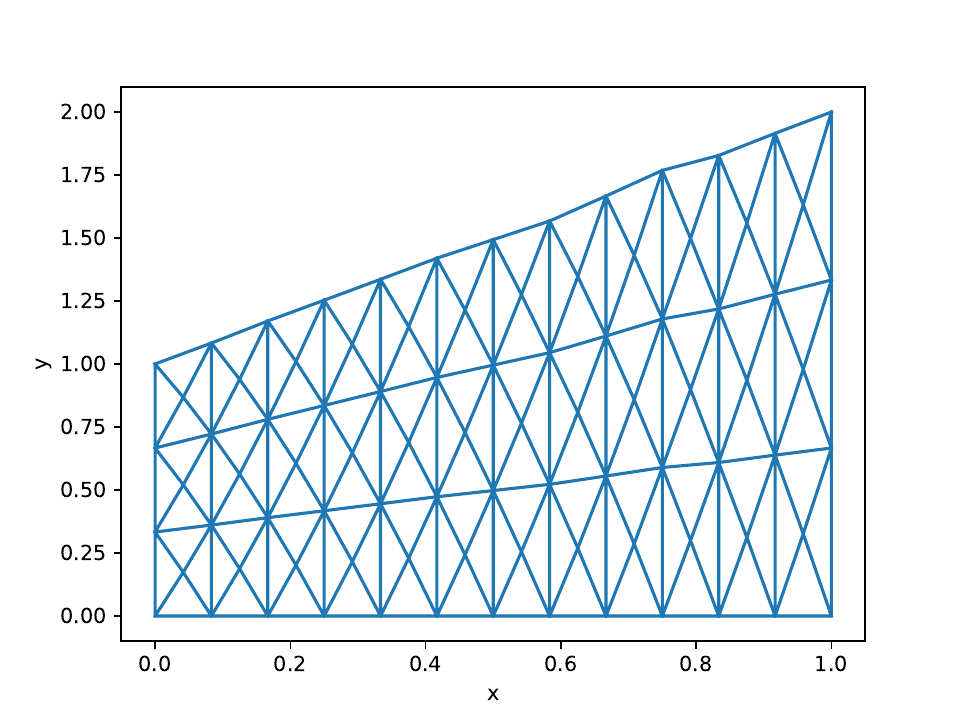}}
  \caption{The initial domain and the change of the domain in three following Newton-like iterations. The free surface is updated vertically.}
\label{domainchange}
\end{figure}

Next we present numerical experiments in 2D.
We start with a straightforward test case for the Dirichlet boundary condition problem and then focus on the submerged triangle problem. The first test case is also a Bernoulli free-boundary problem simplified from the submerged triangle problem, with a Dirichlet condition on both the fixed and free boundary. The submerged triangle problem is the problem to which we are mainly interested in applying this shape-Newton scheme. We will use the algorithm in Table \ref{de} such that the displacement of the free boundary is updated vertically. 

We use a finite element method, based on the weak form of the linearized system in Table~\ref{domainchange}. That is, we seek $(\delta\phi,\delta \eta) \in \hat{V}_h \times \hat{W}_h$ such that
\begin{subequations}
\label{eq:discreteWeakForm}
\begin{multline}
    \int_{\hat{\Omega}}\nabla\delta\phi\cdot\nabla v d \Omega + \int_{\hat{\Omega}\setminus\hat{\Gamma}}\omega\delta\phi vd\Gamma- \int_{\hat{\Gamma}}\text{div}_{\Gamma}\left(\delta\boldsymbol\theta\cdot\boldsymbol n\nabla_{\Gamma}\hat{\phi}\right)v d \Gamma - \int_{\hat{\Gamma}}fv\delta\boldsymbol\theta\cdot\boldsymbol n d \Gamma 
\\ = -\mathcal{R}_1\left(\left(\hat{\boldsymbol\theta}, \hat{\phi}\right); v\right), \quad \forall v\in \hat{V}_h,
\end{multline}
\begin{multline}
    \int_{\hat{\Gamma}}2a\nabla_{\Gamma}\hat{\phi}\cdot\nabla_{\Gamma}\delta\phi w d \Gamma + \int_{\hat{\Gamma}}\left(2\kappa\left|\nabla\hat{\phi}\right|^2+bn_N\right)w\delta\boldsymbol\theta\cdot\boldsymbol n d \Gamma 
\\
= -\mathcal{R}_2\left(\left(\hat{\boldsymbol\theta}, \hat{\phi}\right); w\right), \quad \forall w\in \hat{W}_h.
\end{multline}
\end{subequations}
where $\hat{V}_h$ and $\hat{W}_h$ are finite element spaces based on a quasi-uniform partition (triangulation) of~$\hat{\Omega}$ into a set of shape-regular simplicial elements~$\hat{\Omega}_h$. In particular, we choose continuous piecewise-\emph{linear} approximations, i.e., $\hat{V}_h= \mathbb{P}^1(\hat{\Omega}_h)$ and $\hat{W}_h = \mathbb{P}^1_{0,\mathrm{in}}(\Gamma_h)$, where the (line) elements in $\Gamma_h$ correspond to the free-boundary edges of the (triangular) elements in~$\hat{\Omega}_h$ adjacent to~$\hat{\Gamma}$. The space $\mathbb{P}^1_{0,\mathrm{in}}(\Gamma_h)$ incorporates the condition $\delta\eta =0$ at the inflow of the free boundary (recall Remark~\ref{remark_compatibility1}).
\par
Notice that in the 2D case, 
$$\nabla_{\Gamma}\hat{\phi} = \frac{d\hat{\phi}}{ds}
\qquad \text{and} \qquad 
\delta\boldsymbol\theta\cdot\boldsymbol n = \frac{\delta\eta}{\sqrt{1+ (\hat{\eta}_x)^2}}$$
on the free surface $\hat{\Gamma}$, where $s$ represents the arc length along the free surface, which allows us to write~\eqref{eq:discreteWeakForm} in terms of~$(\delta \phi,\delta \eta)$. 
\par
\begin{remark}[Mesh deformation]
Given a new $\delta\eta$, the free boundary $\hat{\Gamma}$ is updated by moving the mesh nodes on $\hat{\Gamma}$ vertically with the distance~$\delta\eta$. The other mesh nodes are then updated accordingly to yield a smoothly deformed mesh. In particular, we update the other mesh nodes simply by moving vertically using a linearly-interpolated fraction of the distance $\delta\eta$ at the same $x$-coordinate. Further implementation details can be found in Section 6.6.2 in~\cite{yiyunphd}.
\end{remark}
\par
\begin{remark}[Discrete curvature]
For simplicity, the curvature $\kappa$ along the free surface is evaluated by a finite difference approximation. An alternative to this approximation of~$\kappa$ is to obtain the linearisation directly from piecewise smooth free boundaries (which we have not pursued in this work); cf.~Remark \ref{remark_piecewise}.
\end{remark}
\begin{remark}[Solvability of the discrete linear system]\label{remark_solve_discrete}
    Under certain conditions of the mesh and data, the solvability of the discrete shape-Newton schemes for Bernoulli boundary conditions has been proven in Appendix \ref{app:solvability_discrete}.
\end{remark}

\subsection{Dirichlet boundary condition}\label{chap3:numerical_dirichlet}
The test case for the free-boundary problem with Dirichlet boundary condition is a Bernoulli free-boundary problem derived from a manufactured solution,
\begin{equation}
\phi = x +y,  \quad  \eta = x +1, \label{Exact_s}
\end{equation}
such that the data can be obtained as
\begin{eqnarray}
f &=& 0,\nonumber\\
g &=& 0, \nonumber\\
h &=& \begin{cases}
2y -1, \quad \text{on }\Gamma_F,\\
x + y, \quad \text{on }\partial\Omega\setminus\Gamma_F,\nonumber
\end{cases}\\
\omega &\to& \infty\nonumber.
\end{eqnarray}

With an initial domain $\Omega_0 = \left\{\left(x, y\right): x\in\left[0, 1\right], y\in\left[0, x^2+1\right]\right\}$, how the domain and the triangulation changes in the first three iterations are shown in Figure \ref{domainchange}. Starting with a parabola, the free boundary is almost a straight line after the third iteration. The source term has been tested in \cite{van2013shape} by choosing a more complicated manufactured solution.

\begin{figure}[t!]
  \centering
    \includegraphics[width=\linewidth]{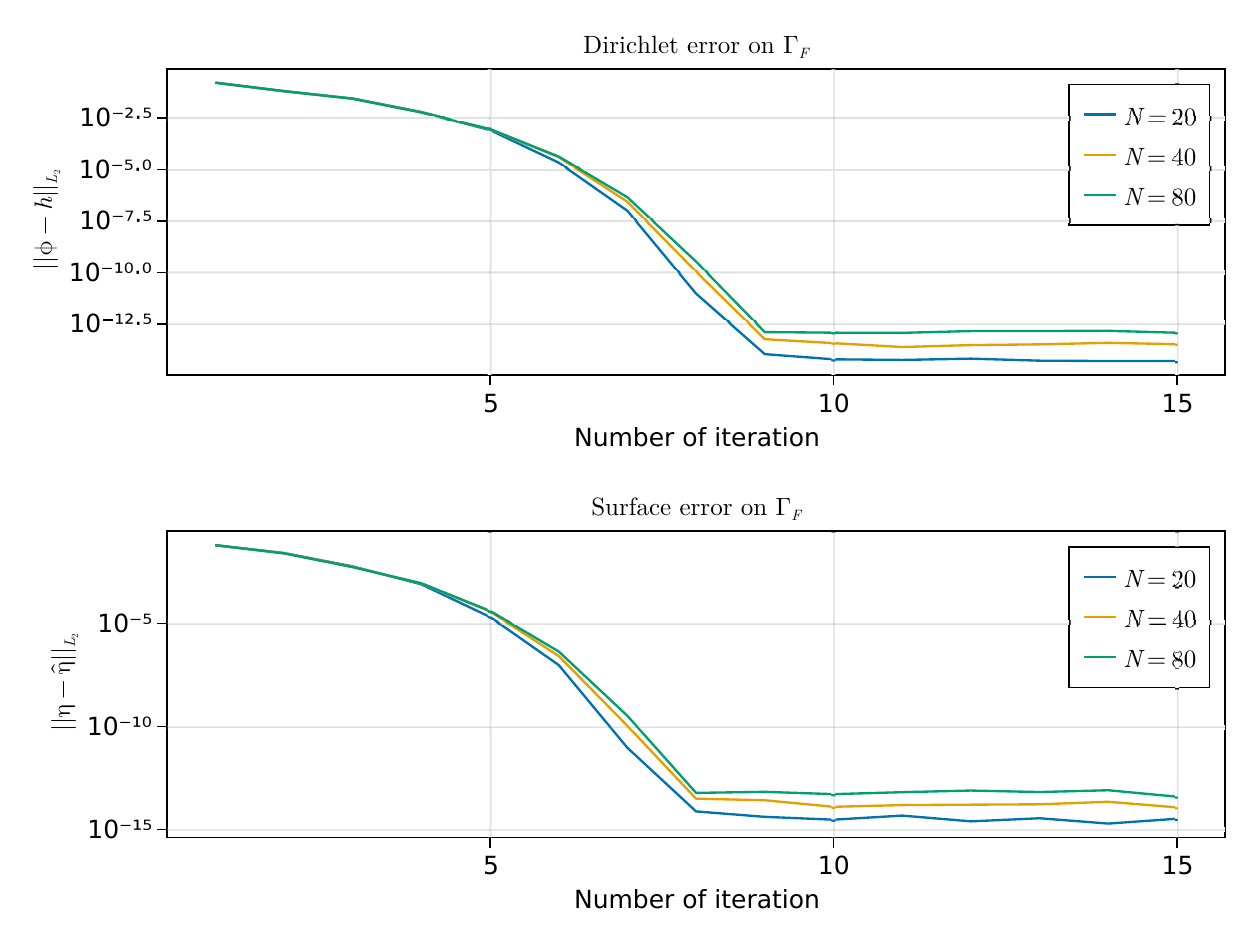}
    \caption{The Dirichlet error $||\phi - h||_{L_2}$ and surface error $|| \eta-\hat{ \eta}||_{L_2}$ on $\Gamma_F$ measured in $L_\infty$-form against the number of iterations. The upper plot shows the Dirichlet error, and the lower shows the surface error. The values of $N+1$ are the number of the nodes along the $x$-axis.}
\label{convergence_D}
\end{figure}

Figure \ref{convergence_D} shows the error between numerical results of $\phi$ and $\eta$ compared with the exact solution \eqref{Exact_s} on the free boundary $\Gamma_F$ with a different number of finite element meshes. The value of $N+1$ represents the number of nodes along the $x$-axis, and the number of nodes along the $y$-axis is $\frac{N}{4}$. Although the error is slightly larger with more nodes, the shape-Newton scheme converges superlinearly.

\subsection{The submerged triangle problem}\label{chap3:numerical_submerged}

\begin{figure}[t!]
\centering
\resizebox{0.75\columnwidth}{!}{%
\begin{tikzpicture}
\draw [very thin] (-4, 0) -- (-1,0);
\draw [very thin] (-1,0) -- (0,1);
\draw [very thin] (0,1) -- (1,0);
\draw [very thin] (1,0) -- (4,0);

\draw[very thin] (-4,2.2) .. controls (0,4) and (0,4) .. (4,2.2);

\draw[very thin] (-4,0) -- (-4,2.2);

\node at (-2.8,2.4)[ inner sep=0pt, scale=0.7] {$\Gamma_F$};
\node at (-4.3,1)[ inner sep=0pt, scale=0.7] {$\Gamma_L$};
\node at (4.3,1)[ inner sep=0pt, scale=0.7] {$\Gamma_R$};
\node at (-3.5,0.2)[ inner sep=0pt, scale=0.7] {$\Omega$};

\draw[ very thin] (4,0) -- (4,2.2);

\draw[dashed] (-1,0) -- (0,0);
\draw[dashed] (0,0) -- (0, 1);

\draw[dashed, |->, thin] (-1,0) -- (0,0);
\node at (-0.5,-0.2)[inner sep=0pt, scale=0.7] {$w_0$};

\node at (3.5,-0.2)[ inner sep=0pt, scale=0.7] {$\Gamma_B$};
\node at (4.2,-0.2)[ inner sep=0pt, scale=0.7] {$(4,0)$};
\node at (-4.2,-0.2)[ inner sep=0pt, scale=0.7] {$(-4,0)$};
\node at (-4.3,2.4)[ inner sep=0pt, scale=0.7] {$(-4,1)$};
\node at (4.3,2.4)[ inner sep=0pt, scale=0.7] {$(4,1)$};

\coordinate (o1) at (-1,0);
\coordinate (a1) at (-0.5,0);
\coordinate (b1) at (0,1);
\coordinate (c1) at (1,0);
\coordinate (d1) at (2,0);
\pic[ draw=black, -, angle eccentricity=1.2, angle radius=0.15cm]
		{angle=a1--o1--b1};
		\node at (-0.5,0.2)[ inner sep=0pt, scale=0.6] {$\alpha$};

\end{tikzpicture}%
}
\caption{The sketch of the domain we used for the second test case. $\alpha$ is denoted as the angle and $w_0$ as the half width of the triangle.}
\label{D_2}
\end{figure}
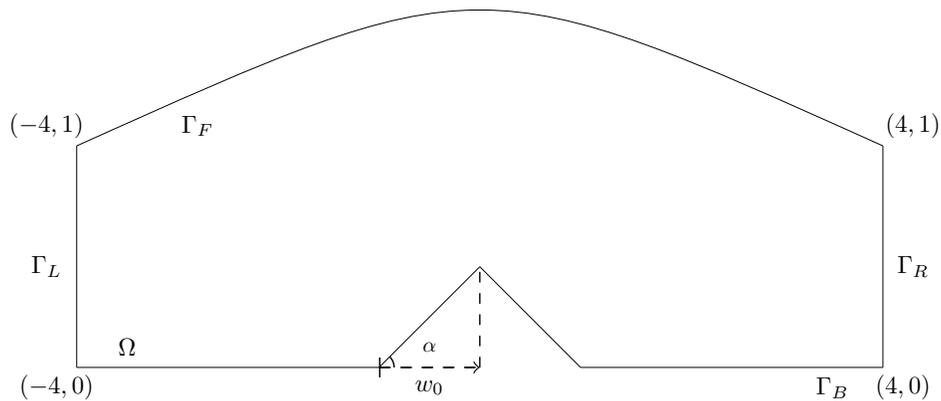
\begin{figure}[t!]
  \centering
    \includegraphics[width=\linewidth]{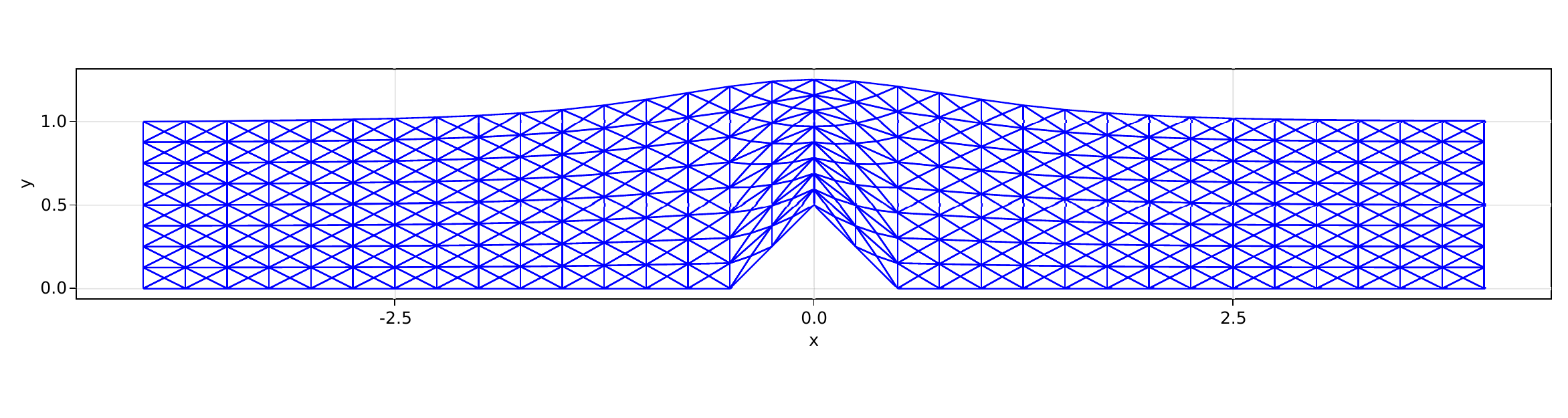}
    \caption{An example of the domain and the triangulation with $\alpha = \frac{\pi}{4}$, $F=2$ and the half width of the triangle $w_0 = 0.5$.}
\label{domain_B}
\end{figure}

The second test case is the submerged triangle problem investigated by Dias and Vanden-Broeck \cite{dias1989open}. A detailed derivation of the governing equations can be found in \cite[Appendix]{yiyunphd}. In this problem, we have a Neumann boundary condition on $\partial\Omega\setminus\Gamma_R$ and a Dirichlet boundary condition on $\Gamma_R$, i.e. $\omega = 0$ on $\partial\Omega\setminus\Gamma_R$ and $\omega\to\infty$ on $\Gamma_R$. The data defining this problem is given as follows:
\begin{eqnarray}
f &=& 0,\nonumber\\
g &=& \begin{cases}
0, \quad \text{on }\partial\Omega\setminus\Gamma_L,\\
-1, \quad \text{on }\Gamma_L,
\end{cases}\nonumber\\
h &=& 0 \quad \text{on } \Gamma_R.\nonumber
\end{eqnarray}
The Bernoulli condition is obtained by giving $a = \frac{1}{2}F^2$, $b =1$ and $c = - \frac{1}{2}F^2-1$ where $F$ is the Froude number. The domain is a rectangle truncated at $|x|=4$ containing an isosceles triangle symmetric about $x=0$ having an angle $\alpha$ and width $2w_0$ at the bottom, as shown in Figure \ref{D_2}. The space is discretised as shown in Figure \ref{domain_B}, where it was uniformly spaced along the $x-$axis and the vertical direction for fixed values of $x$. Then the algorithm in Table \ref{de} can be applied to solve for the pair $\left(\delta\eta, \delta \phi\right)$, and the free boundary can be updated vertically with $\delta\eta$.

Dias and Vanden-Broeck \cite{dias1989open} found that the solutions to the submerged problem have two types: One is supercritical flow both upstream and downstream, and the other is supercritical (or subcritical) upstream and subcritical (or supercritical) downstream flow. Our numerical solutions are the first type, and we can compare them with the results in \cite{dias1989open}.

\subsubsection{Convergence rate of Shape-Newton method}
\begin{figure}[t!]
  \centering
    \includegraphics[width=\linewidth]{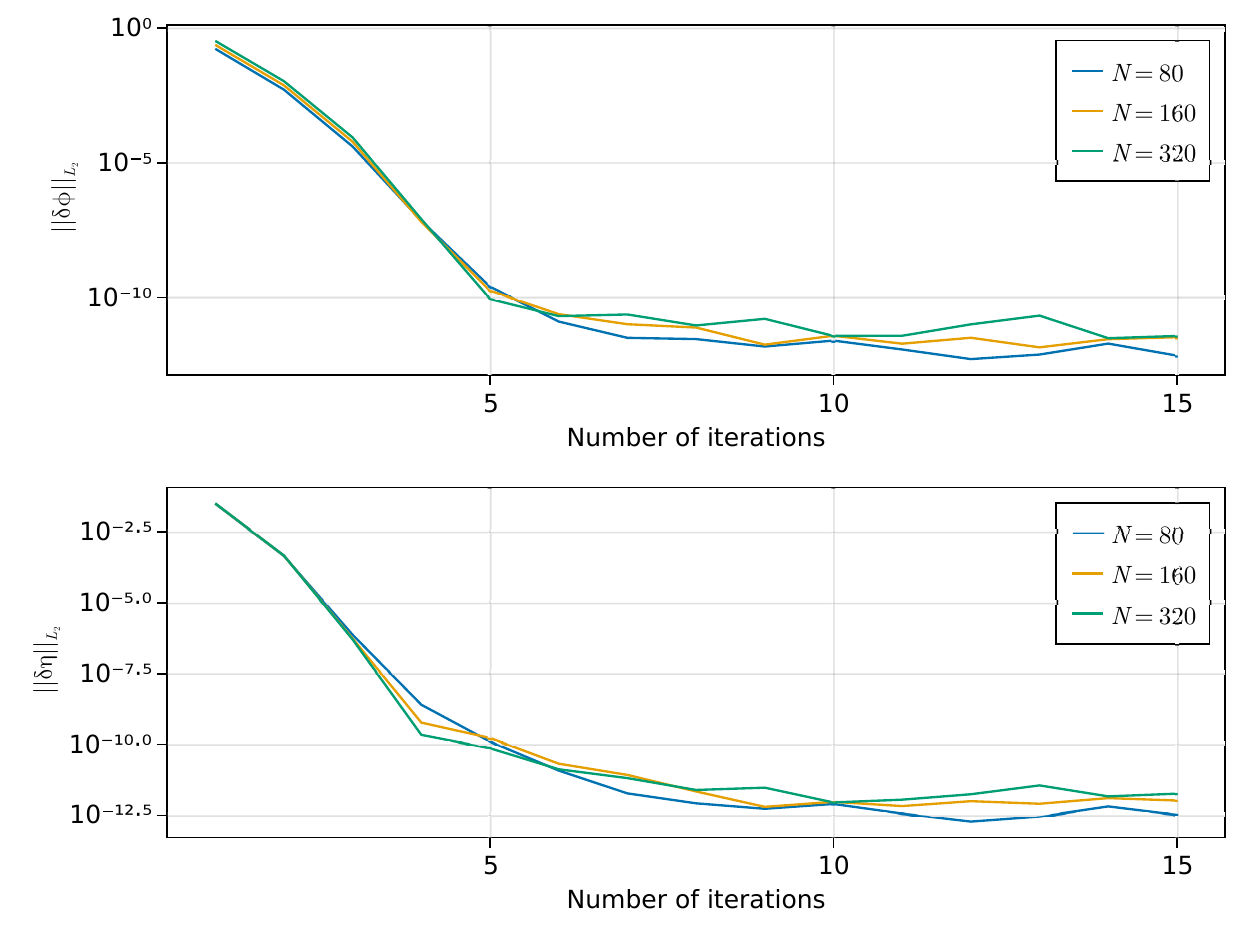}
    \caption{The size of $||\delta\phi||_{L_2}$ and $|| \delta \eta||_{L_2}$ on $\Gamma_F$ against the number of iterations with $\alpha = \frac{\pi}{8}$, $w_0 = 0.3$ and $F=3$. The values of $N+1$ are the number of the nodes along the $x$-axis. }
\label{Err_B}
\end{figure}

The rate of convergence is shown in Figure \ref{Err_B}, where we show $||\delta\phi||_{L_2}$ and $||\delta \eta||_{L_2}$ against the number of iterations for $\alpha = \frac{\pi}{8}$, $w_0 = 0.3$ and $F=3$. These show superlinear convergence. This figure also shows the comparison for different mesh densities.

\subsubsection{Robustness of the Shape-Newton scheme}
\begin{figure}[t!]
  \centering
  \subfloat[The final domain for $\alpha=\frac{\pi}{16}$, $w_0=0.5$, and $F=2$.]{\includegraphics[width=\linewidth]{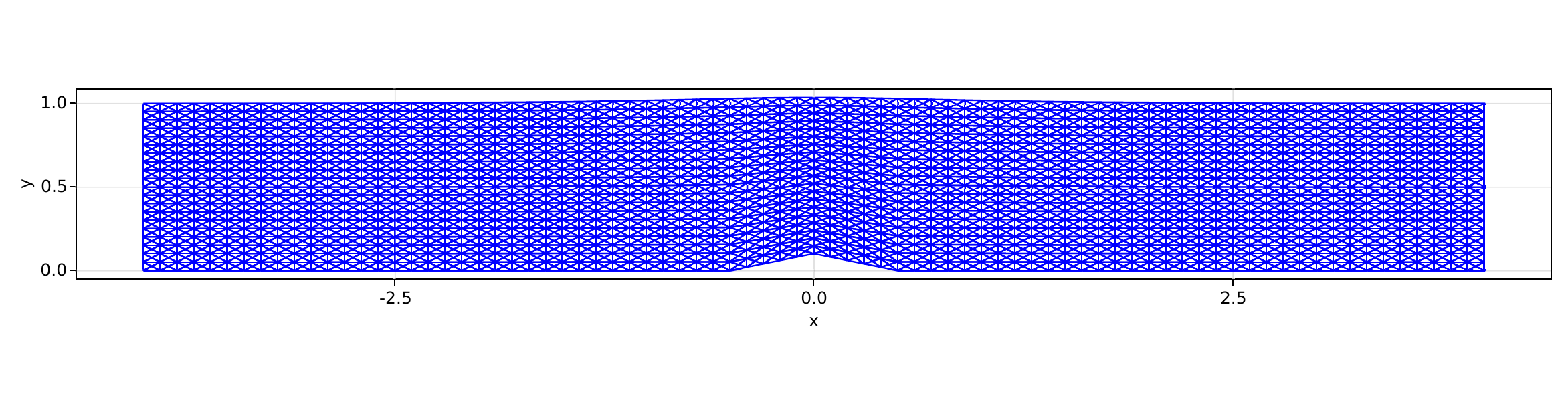}}\hfill
  
  \subfloat[The final domain for $\alpha=\frac{\pi}{8}$, $w_0=0.5$, and $F=1.4$.]{\includegraphics[width=\linewidth]{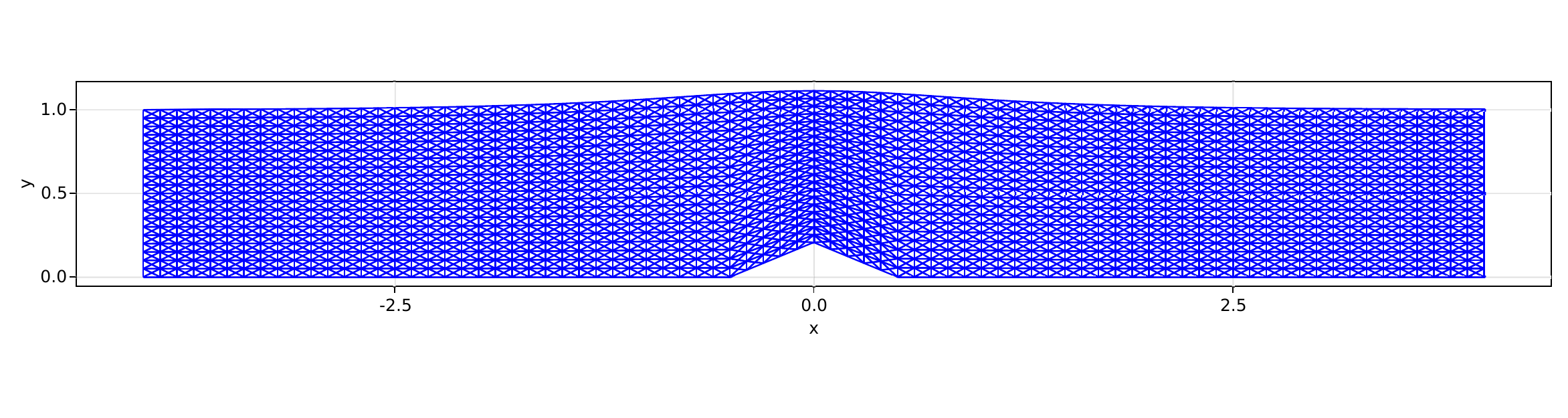}}\hfill
  
    \subfloat[The final domain for $\alpha=\frac{\pi}{8}$, $w_0=0.5$, and $F=2$.]{\includegraphics[width=\linewidth]{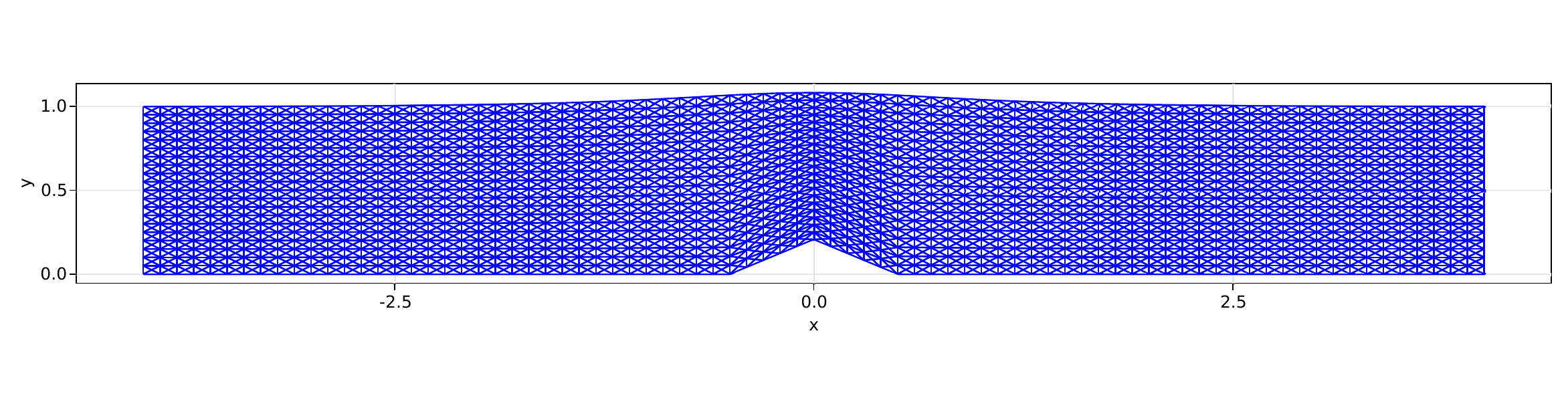}}\hfill
    
    \subfloat[The final domain for $\alpha=\frac{\pi}{8}$, $w_0=1$, and $F=2$.]{\includegraphics[width=\linewidth]{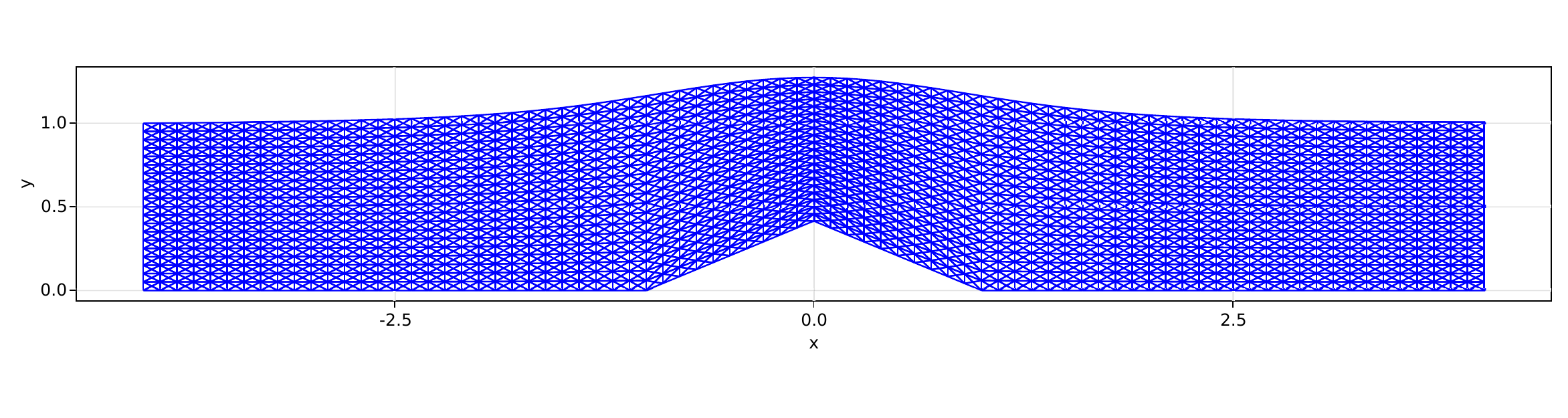}}
  \caption{The final domains for various $\alpha$, $w_0$ and $F$, where their free boundaries are the numerical solutions. }
\label{solu}
\end{figure}

\begin{figure}[t!]
  \centering
  \subfloat[The maximum value $ y$ on the free boundary at $x=0$ against $F$ with $w_0=0.1$ for different values of $\alpha$.]{\label{y0_d1}\includegraphics[width=0.5\linewidth]{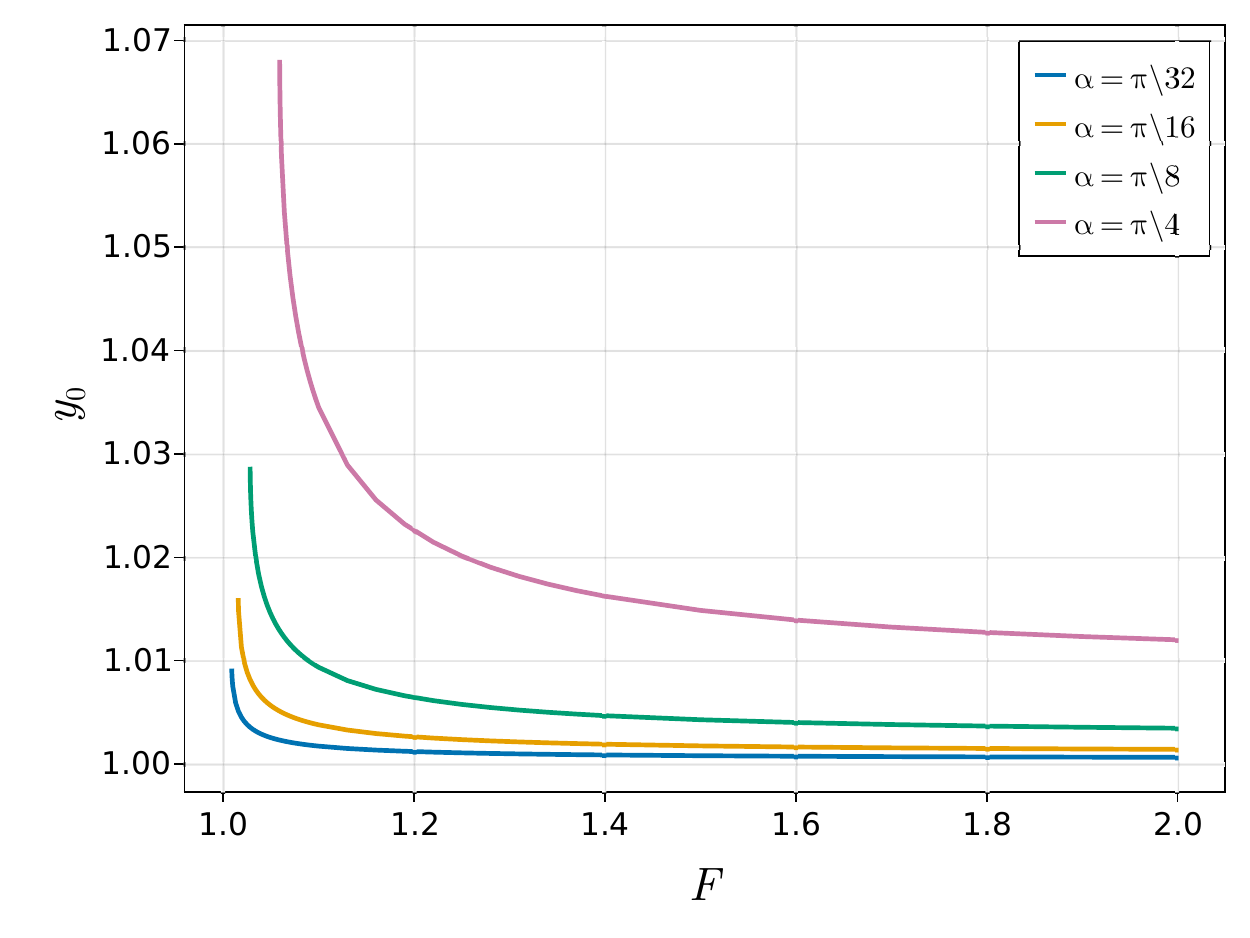}}\hfill
    \subfloat[The maximum value $ y$ on the free boundary at $x=0$ against $F$ with $w_0=0.3$ for different values of $\alpha$.]{\label{y0_d3}\includegraphics[width=0.5\linewidth]{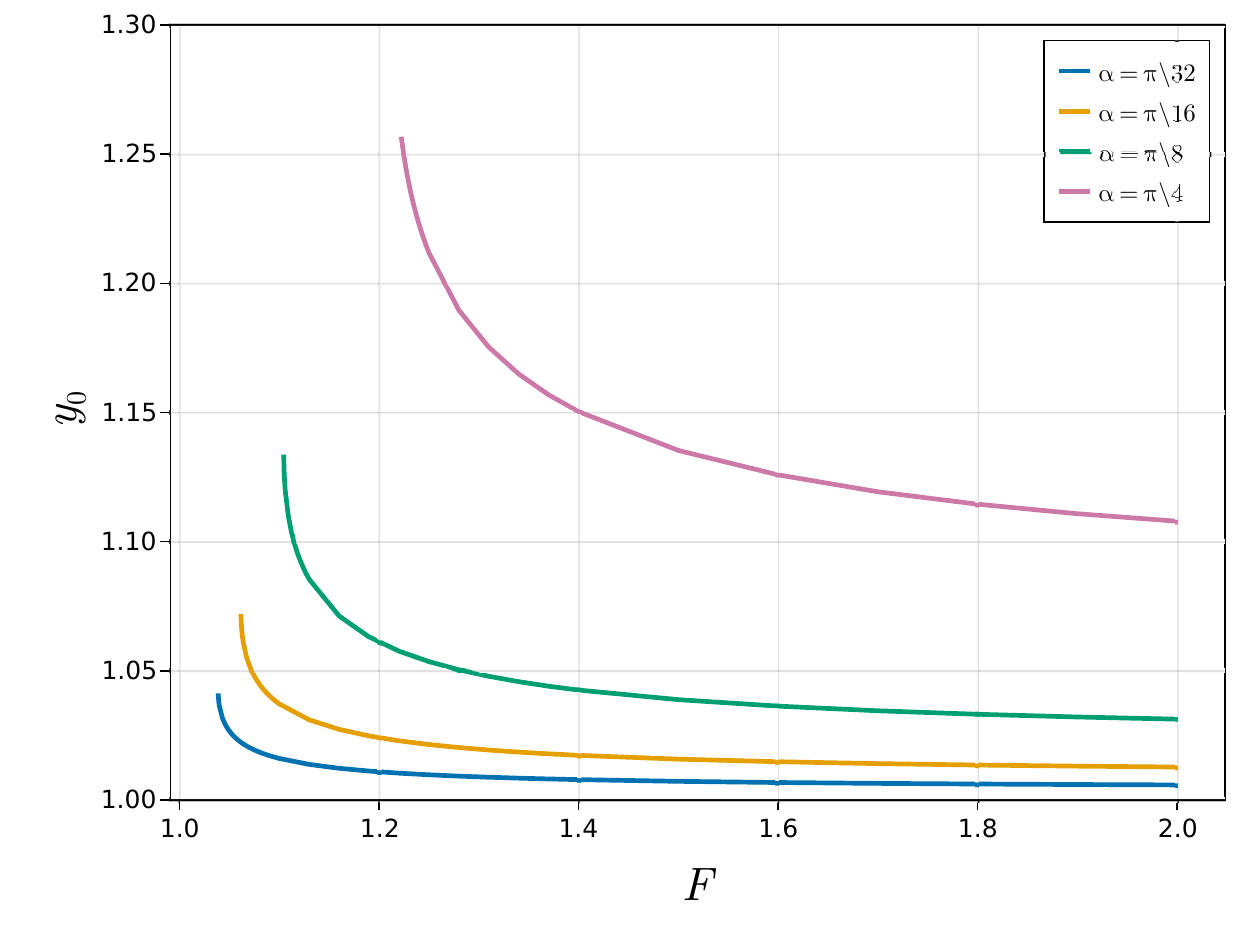}}
    
    \subfloat[The maximum value $ y$ on the free boundary at $x=0$ against $F$ with $w_0=0.5$ for different values of $\alpha$.]{\label{y0_d5}\includegraphics[width=0.5\linewidth]{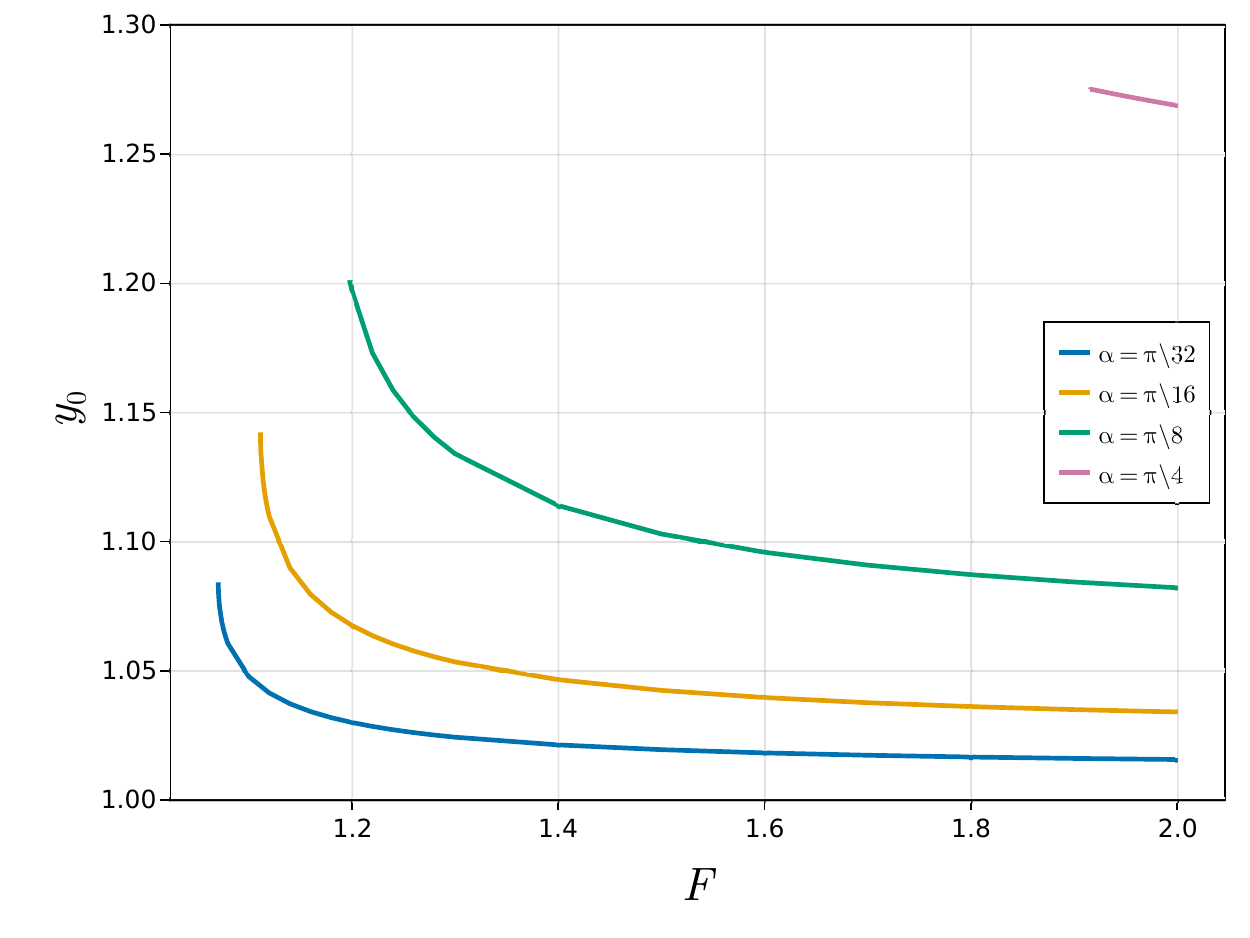}}\hfill
    \subfloat[The maximum value $ y$ on the free boundary at $x=0$ against $F$ with $w_0=0.7$ for different values of $\alpha$.]{\label{y0_d7}\includegraphics[width=0.5\linewidth]{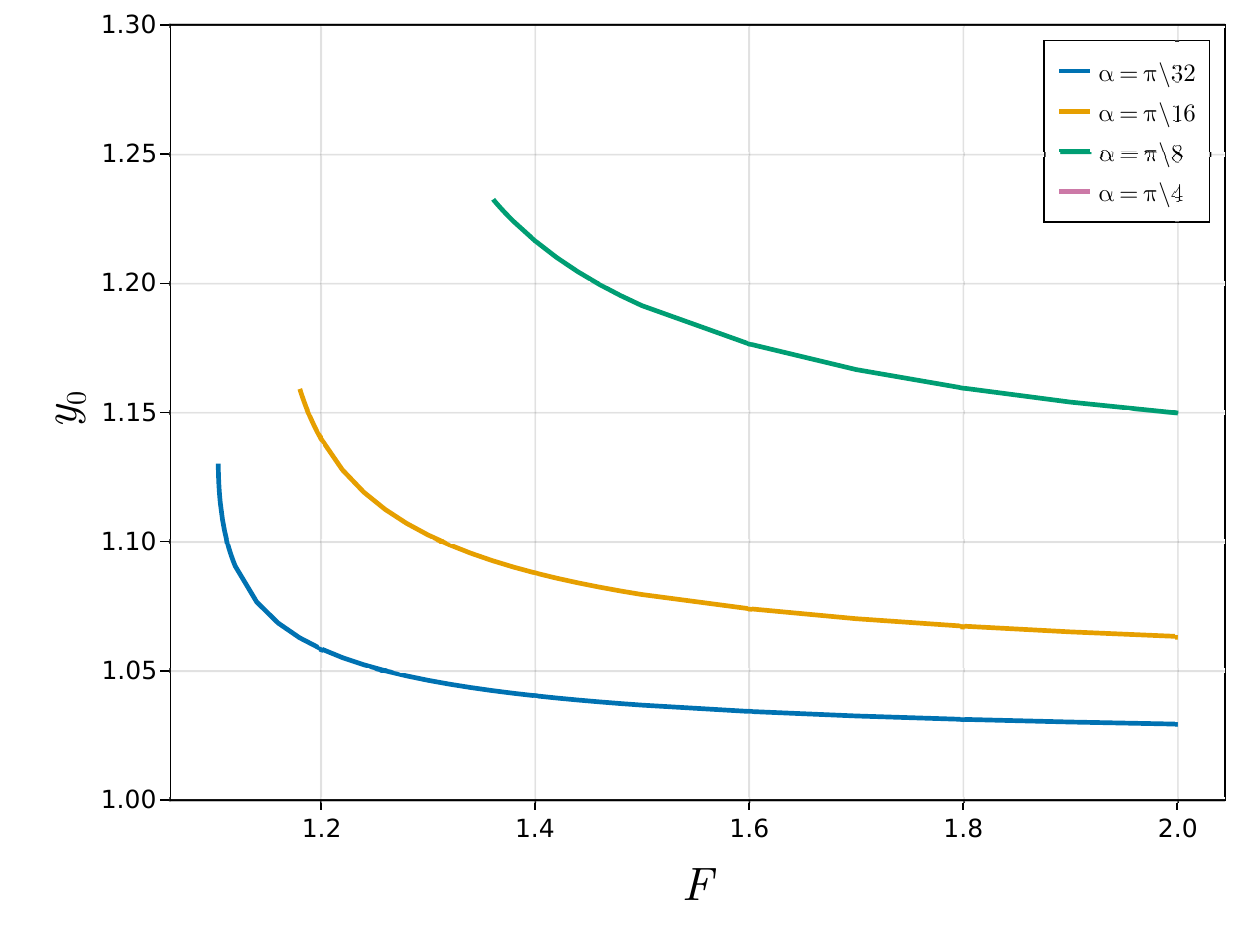}}

  \caption{The maximum value $ y_0$ on the free boundary at $x=0$ against $F$ for different values of $\alpha$ and $w_0$. }
\label{crit}
\end{figure}

Some converged grids of the whole region are shown in Figure \ref{solu}. We noticed that $\eta(x)$ has a maximum value $y_0$ at $x=0$ on the free boundary, and the value of $y_0$ changes with the values of $\alpha$, $w_0$ and $F$. Figure \ref{crit} shows the value of $ y_0$ against the Froude number $F$ for various values of $\alpha$. We can observe from Figure \ref{crit} that $ y_0$ will decrease when the Froude number $F$ becomes larger for the fixed width of the triangle. In addition, for fixed values of $F$ and angle $\alpha$, $y_0$ will also decrease with the width of the triangle. This agrees with the results presented by Dias and Vanden-Broeck in \cite{dias1989open}, who solved this problem for fixed $\alpha=\frac{\pi}{4}$. To improve convergence behaviour, we explored using a continuation technique in the Froude number F. However, as seen in Figure \ref{y0_d5} and Figure \ref{y0_d7}, for larger triangles (larger $\alpha$), convergence generally becomes more difficult, as those test cases are closer to critical situations beyond which there is no solution. See \cite{dias1989open} for a detailed study on critical values.

We also found that the solutions are challenging for larger angle $\alpha$ for fixed width. The possible reason is that with a higher triangle height, the flow can approach its limiting configuration as a thin layer over the edge of the triangle with a stagnation point, hence may require local mesh refinement.

\section{Conclusion}\label{chap:3_conclusion}
We derived a shape-Newton method to solve generic \break free-boundary problems with the nonlinear Bernoulli boundary condition. The linearised system is obtained from applying the Hadamard formula for shape derivatives to a suitable weak form of the free boundary problem. After linearisation and neglecting higher-order terms, one obtains a linear boundary-valued problem to be solved at each iteration.

The shape linearisation of the nonlinear Bernoulli equation is a key result in our work. In its derivation, many terms can be neglected (as a higher-order correction) due to the homogeneous Neumann boundary condition. After some calculations, we find that the result involves the normal derivative of the velocity squared, i.e. $\partial_{\boldsymbol n}\left|\nabla\phi\right|^2$. This can be equivalently computed as $\left(\nabla\hat{\phi}\right)^T\cdot\left[\nabla_{\Gamma}\boldsymbol n\right]\cdot\nabla\hat{\phi}$, see Section~\ref{linearisation_bernoulli}. 

The linearised system essentially corresponds to a boundary value problem for the Laplacian with a generalized Robin boundary condition involving a surface Laplacian (Laplace--Beltrami), which in turn depends on the curvature. Another key result in our work is a study of the solvability of this linearised system. Under certain conditions on the data, one can guarantee the existence of a unique solution (details in Appendix~\ref{app:solvability_bernoulli}).
We applied our method to compute the flow over a submerged triangle for a range of Froude numbers and triangle shapes, and obtained consistent results with the earlier literature \cite{dias1989open}. Moreover, the numerical test revealed that the shape-Newton method converges superlinearly. A theoretical explanation of this behaviour remains an open problem.

\appendix
\section{Consistency of \texorpdfstring{\eqref{dn}}{dn} with \texorpdfstring{\eqref{general_deriviation}}{de}}\label{Appendix: calculation}
In this appendix, we will show the detail about consistency between the result \eqref{general_deriviation} for the $N$-dimensional case and the result \eqref{dn} for the two-dimensional case. 

In the two-dimensional case, we have the unit normal vector  $\boldsymbol n(x) =\left(n_x,n_y\right)=\frac{1}{\sqrt{1+\hat{\eta}^2_x}}\left(-\hat{\eta}_x,1\right)$. Thus, 
\begin{eqnarray}
    \nabla\boldsymbol n &=& \begin{pmatrix}
        \partial_x n_x & \partial_x n_y\\
        \partial_y n_x & \partial_y n_y
    \end{pmatrix} = \begin{pmatrix}
        -\kappa & -\hat{\eta}_x\kappa\\
        0 &0
    \end{pmatrix},\label{appedix: grad_n}\\
    \left(\boldsymbol n\cdot \nabla\right)\boldsymbol n &=& \left(n_x\partial_x + n_y\partial_y\right)\boldsymbol n = n_x\partial_x\left(\boldsymbol n\right) = \begin{pmatrix}
        -\kappa n_x\\
        -\hat{\eta}_xn_x\kappa
    \end{pmatrix},\label{appendix:gradn_n}
\end{eqnarray}
where $\kappa = \partial_x\left(\frac{\hat{\eta}_x}{\sqrt{1+\hat{\eta}^2_x}}\right) = -\partial_xn_x$ represents the curvature. Hence,
\begin{eqnarray}
    \partial_{\boldsymbol n}\left(\boldsymbol n\right) \boldsymbol n^T =\begin{pmatrix}
        -\kappa n_x \\
        -\kappa\hat{\eta}_x n_x
    \end{pmatrix}\begin{pmatrix}
        n_x & n_y
    \end{pmatrix} = \begin{pmatrix}
        -\kappa n_x^2 & -\kappa n_xn_y\\
        -\kappa\hat{\eta}_xn_x^2 & -\kappa\hat{\eta}_xn_xn_y
    \end{pmatrix}\label{appendix:gradn_nn}
\end{eqnarray}

By substituting \eqref{appedix: grad_n} and \eqref{appendix:gradn_nn} into the definition of $\nabla_{\Gamma}\left(\cdot\right)$ \eqref{tangential_gd}, we obtain
\begin{eqnarray}
    \nabla_{\Gamma}\boldsymbol n &=& \nabla \boldsymbol n - \partial_{\boldsymbol n}\boldsymbol n\boldsymbol n^T\nonumber\\
    &=&\begin{pmatrix}
        -\kappa & -\hat{\eta}_x\kappa\\
        0 & 0
    \end{pmatrix}- \begin{pmatrix}
        -\kappa n_x^2 & -\kappa n_xn_y\\
        -\kappa\hat{\eta}_xn_x^2 &-\kappa\hat{\eta}_xn_xn_y
    \end{pmatrix}\nonumber\\
    &=& \kappa\begin{pmatrix}
        -1+n_x^2 & -\hat{\eta}_x+n_xn_y\\
        \hat{\eta}_xn_x^2 & \hat{\eta}_xn_xn_y
    \end{pmatrix}.\label{appendix: gardgamma_n}
\end{eqnarray}
Now by using \eqref{appendix: gardgamma_n} and the Neumann boundary condition \eqref{FB_N}, we have
\begin{alignat}{2}
\left(\nabla_\Gamma\hat{\phi}\right)^T\cdot\left[\nabla_{\Gamma}\boldsymbol n\right]\cdot\nabla_{\Gamma}\hat{\phi} &= \left(\nabla\hat{\phi}\right)^T\cdot\left[\nabla_{\Gamma}\boldsymbol n\right]\cdot\nabla\hat{\phi}\nonumber\\
&= \kappa\begin{pmatrix}
    \hat{\phi}_x & \hat{\phi_y}
\end{pmatrix}\begin{pmatrix}
        -1+n_x^2 & -\hat{\eta}_x+n_xn_y\\
        \hat{\eta}_xn_x^2 & \hat{\eta}_xn_xn_y
    \end{pmatrix}\begin{pmatrix}
        \hat{\phi}_x\\
        \hat{\phi}_y
    \end{pmatrix}\nonumber\\
    &= \kappa\begin{pmatrix}
        \left(-1+n_x^2\right)\hat{\phi}_x + \hat{\eta}_xn_x^2\hat{\phi}_y\\
        \left(-\hat{\eta}_x+n_xn_y\right)\hat{\phi}_x + \hat{\eta}_xn_xn_y\hat{\phi}_y
    \end{pmatrix}\begin{pmatrix}
        \hat{\phi}_x\\
        \hat{\phi}_y
    \end{pmatrix}\nonumber\\
    &= \kappa\bigl[\left(-1+n_x^2\right)\hat{\phi}_x^2+\hat{\eta}_xn_x^2\hat{\phi}_x\hat{\phi}_y + \left(-\hat{\eta}_x+n_xn_y\right)\hat{\phi}_x\hat{\phi}_y \nonumber\\
    &\qquad+ \hat{\eta}_xn_xn_y\hat{\phi}_y^2\bigr]\nonumber\\
    &= \kappa\bigl[\left(-1+n_x^2\right)\hat{\phi}_x^2+\hat{\eta}^2_xn_x^2\hat{\phi}_x^2 + \left(-\hat{\eta}_x+n_xn_y\right)\hat{\eta}_x\hat{\phi}_x^2 \nonumber\\
    &\qquad+ \hat{\eta}_x^3n_xn_y\hat{\phi}_x^2\bigr]\tag{by \eqref{FB_N}}\\
    &=\kappa\hat{\phi}^2_x\left(-1+n_x^2+\hat{\eta}_x^2n_x^2-\hat{\eta}_x^2+\hat{\eta}_xn_xn_y+\hat{\eta}^3_xn_xn_y\right)\nonumber\\
    &= \kappa\hat{\phi}^2_x\left(-1+n_x^2+\hat{\eta}_x^2n_x^2-\hat{\eta}^2_x-n_x^2-\hat{\eta}_x^2n_x^2\right)\tag{by $n_x=-\hat{\eta}_xn_y$}\\
    &= -\left(1+\hat{\eta}_x^2\right)\kappa\hat{\phi}_x^2\nonumber\\
     &=\frac{2}{\sqrt{1+\hat{ \eta}_x^2}}\hat{ \eta}_{xx}\left(\hat{\phi}_x\right)^2\nonumber\\
    &= -\kappa\left|\nabla\hat{\phi}\right|^2,\tag{by \eqref{FB_N}}
\end{alignat}
    Hence, \eqref{general_deriviation} in the two-dimensional case equals \eqref{dn}.

\section{Solvability of the shape-linearized system for the Dirichlet boundary condition}\label{app:solvability_dirichlet}
In this Appendix we show that, under certain conditions of the data, the shape-linearized system~\eqref{SF_laplace_D}--\eqref{CS_D} for the free-boundary problem with Dirichlet boundary condition (i.e.,~\eqref{FBD_L}--\eqref{FBD_LBC}), has a unique solution.
    \par
From the Dirichlet boundary condition \eqref{CS_D}, we have
\begin{equation}
    \delta\boldsymbol{\theta}\cdot\boldsymbol{n} = \frac{\delta\phi+\hat{\phi}-h}{\partial_{\boldsymbol{n}}h}, \label{dtheta_d1}
\end{equation}
provided~$\partial_{\boldsymbol{n}}h\neq0$. Note that for the case when $h$ is a constant and $\partial_{\boldsymbol{n}}h=0$, this problem has been shown to have a unique solution; see, e.g.~\cite{van2013shape}.

Substituting~\eqref{dtheta_d1} into \eqref{CS_ND}, the system becomes one-way coupled, i.e., a \hfill\break boundary-value problem for~$\delta \phi$, and subsequently an equation for~$\delta \boldsymbol{\theta}\cdot \boldsymbol{n}$:
\begin{subequations}
\begin{gather}
\nabla^2\delta\phi = -\nabla^2\hat{\phi}-f\quad \text{in }\hat{\Omega},
\label{DCSD_L}
\\
\quad\partial_n\delta\phi -\text{div}_{\Gamma}\left(\frac{\delta\phi}{\partial_{\boldsymbol{n}}h}\nabla_{\Gamma}h\right) - f\frac{\delta\phi}{\partial_{\boldsymbol{n}}h} = \text{div}_{\Gamma}\left(\frac{\hat{\phi}-h}{\partial_{\boldsymbol{n}}h}\nabla_{\Gamma}h\right)+f\frac{\hat{\phi}-h}{\partial_{\boldsymbol{n}}h}
-\partial_{\boldsymbol{n}}\hat{\phi} , \quad \text{on }\hat{\Gamma}, 
\label{DCSD_N}
\\
\partial_{\boldsymbol n}\delta\phi + \omega\delta\phi = g+\omega h-\left(\partial_{\boldsymbol n}\hat{\phi} + \omega\hat{\phi}\right) \quad \text{on }\partial\hat{\Omega}\setminus\hat{\Gamma},
\label{DCSD_R}
\end{gather}
\end{subequations}
The boundary-value problem for~$\delta \phi$ has essentially a generalized Robin boundary condition involving an oblique derivative ($\partial_n \delta \phi - \frac{1}{\partial_{\boldsymbol{n}} h} \nabla_\Gamma h \cdot  \nabla_\Gamma \delta \phi$) on~$\hat{\Gamma}$. To guarantee existence of a unique weak solution to this boundary-value problem, we use the Lax-Milgram theorem and establish coercivity of a suitable bilinear form. We assume~$\hat{\Omega}$ is a bounded Lipschitz-continuous domain, and $\hat{\phi}$ is sufficiently smooth to ensure continuity of the bilinear form and linear form of the weak formulation.
\par
A weak form of~\eqref{DCSD_L}--\eqref{DCSD_R} seeks $\delta\phi \in H^1(\hat{\Omega})$ such that
\begin{equation}
    a(\delta\phi,\psi)=l(\psi),
    \qquad \forall \psi\in H^1(\hat{\Omega})\,,
\end{equation}
where
\begin{alignat}{2}    
\label{a_d_uv}
    a(\delta\phi,\psi) &:= \int_{\hat{\Omega}}\nabla\delta\phi\cdot\nabla \psi d\Omega 
+ \int_{\hat{\Omega}\setminus\hat{\Gamma}}\omega\delta\phi \psi d\Gamma
\\ 
&\qquad 
+\int_{\hat{\Gamma}}\frac{\delta\phi}{\partial_{\boldsymbol{n}}h}\nabla_{\Gamma}h \cdot \nabla_\Gamma \psi d\Gamma 
- \int_{\hat{\Gamma}}f\frac{\delta\phi}{\partial_{\boldsymbol{n}}h} \psi d\Gamma,
    \nonumber
    \\
    \label{linear_d_right}
     l(\psi) &:= 
-
\int_{\hat{\Omega}}\nabla\hat{\phi}\cdot\nabla \psi d\Omega + \int_{\hat{\Omega}}f\psi d\Omega
+
\int_{\partial\hat{\Omega}\setminus\hat{\Gamma}}(g+\omega h-\omega\hat{\phi})\psi d\Gamma 
\\
& \qquad 
-\int_{\hat{\Gamma}}\frac{\hat{\phi}-h}{\partial_{\boldsymbol{n}}h}\nabla_{\Gamma}h \cdot \nabla_{\Gamma} \psi
d\Gamma
+\int_{\hat{\Gamma}}f\frac{\hat{\phi}-h}{\partial_{\boldsymbol{n}}h} \psi d\Gamma.
     \nonumber
\end{alignat}
To study the coercivity of~$a(\cdot,\cdot)$, note that
\begin{alignat}{2}
    a(\psi,\psi) &= \int_{\hat{\Omega}} |\nabla\psi|^2 d\Omega + \int_{\partial\hat{\Omega}\setminus\hat{\Gamma}}\omega \psi^2 d\Gamma
    + \int_{\hat{\Gamma}}\frac{\psi}{\partial_{\boldsymbol{n}}h}\nabla_{\Gamma}h\cdot\nabla_{\Gamma}\psi d\Gamma 
    -
    \int_{\hat{\Gamma}}\frac{f}{\partial_{\boldsymbol{n}}h} \psi^2 d\Gamma,
\end{alignat}
hence, if the last two terms are nonnegative, coercivity holds when $\omega >0$ (by a Poincar\'e--Steklov inequality; see, e.g.,~\cite[Eq.~(31.23)]{ern2021finite}). We note that the penultimate term can be written as
\begin{alignat*}{2}
 \int_{\hat{\Gamma}}\frac{1}{2\partial_{\boldsymbol{n}}h}\nabla_{\Gamma}h\cdot\nabla_{\Gamma}(\psi^2)d\Gamma\nonumber
    &=-\int_{\hat{\Gamma}}\text{div}_{\Gamma}\left(\frac{1}{2\partial_{\boldsymbol{n}}h}\nabla_{\Gamma}h\right)\psi^2d\Gamma 
    + \int_{\partial\hat{\Gamma}}\frac{\nabla_{\Gamma}h\cdot\boldsymbol{\tau}}{2\partial_{\boldsymbol{n}}h}\psi^2 d s \,.
    \end{alignat*}
Therefore, sufficient conditions that guarantee coercivity are:
$$\text{div}_{\Gamma}\left(\frac{1}{2\partial_{\boldsymbol{n}}h}\nabla_{\Gamma}h\right)\le 0, \quad \text{and } \frac{f}{\partial_{\boldsymbol{n}}h}\le 0,$$ 
and a closed free boundary $\hat{\Gamma}$ (hence $\partial\hat{\Gamma}=\emptyset$) or $\frac{\nabla_{\Gamma}h\cdot\boldsymbol{\tau}}{2\partial_{\boldsymbol{n}}h}\geq0$.
\par
These conditions then guarantee that $\delta \phi \in H^1(\hat{\Omega})$ while $\delta \boldsymbol{\theta} \cdot \boldsymbol{n} $ follows from~\eqref{dtheta_d1}. Generally, one expects additional regularity for $\delta \phi$ beyond~$H^1(\hat{\Omega})$, so that $\delta \boldsymbol{\theta}$ inherits this regularity and becomes a Lipschitz-continuous vector field on~$\Gamma_0$. Such regularity study is outside the scope of this work.

\section{Solvability of the shape-linearized system for the \break Bernoulli boundary condition}\label{app:solvability_bernoulli}
In this Appendix we show that, under certain conditions of the data, the shape-linearized system~\eqref{linSys:Bern} for the free-boundary problem with Bernoulli boundary condition (i.e.,~\eqref{FB_L}--\eqref{FB_Fixed}), has a unique solution.
    \par
The linearized Bernoulli condition \eqref{CS_B} can be rearranged to:
\begin{equation}
    \label{linBC} \delta\boldsymbol{\theta}\cdot\boldsymbol{n} = \frac{2a\nabla_{\Gamma}\hat{\phi}\cdot\nabla_{\Gamma}\delta\phi+C_2}{C_1}
\end{equation}
where $C_1 = -2\left(\nabla_{\Gamma}\hat{\phi}\right)^T\cdot\left[\nabla_{\Gamma}\boldsymbol n\right]\cdot\nabla_{\Gamma}\hat{\phi}+bn_N$ and $C_2 = a\left|\nabla\hat{\phi}\right|^2+b\hat{x}_N+c$, provided that $C_1\neq 0$.
Let $\boldsymbol{b}_{\Gamma} = \frac{2a\nabla_{\Gamma}\hat{\phi}}{C_1}$ and $C = \frac{C_2}{C_1}$ for notation convenience, \eqref{linBC} can be rewritten as
\begin{equation}
    \delta\boldsymbol{\theta}\cdot\boldsymbol{n} = \boldsymbol{b}_{\Gamma}\cdot\nabla_{\Gamma}\delta\phi + C.\label{dthetan}
\end{equation}

Similar to the approach in Appendix~\ref{app:solvability_dirichlet}, by substituting \eqref{dthetan} into~\eqref{CS_N}, we obtain a boundary-value problem for~$\delta\phi$:
\begin{subequations}
\begin{alignat}{2}
-\nabla^2\delta\phi = \nabla^2\hat{\phi}+f\quad \text{in }\hat{\Omega},\label{DCS_L}\\
\qquad\partial_n\delta\phi+\partial_{\boldsymbol{n}}\hat{\phi} -\text{div}_{\Gamma}\left(\boldsymbol{b}_{\Gamma}\cdot\nabla_{\Gamma}\delta\phi\nabla_{\Gamma}\hat{\phi}\right)- f\boldsymbol{b}_{\Gamma}\cdot\nabla_{\Gamma}\delta\phi = -\text{div}_{\Gamma}\left(C\nabla_{\Gamma}\hat{\phi}\right)-Cf \quad \\
\text{on }\hat{\Gamma},\nonumber\\
\partial_{\boldsymbol n}\delta\phi + \omega\delta\phi = g+\omega h-\left(\partial_{\boldsymbol n}\hat{\phi} + \omega\hat{\phi}\right) \quad \text{on }\partial\hat{\Omega}\setminus\hat{\Gamma},\label{DCS_R}
\end{alignat}
\end{subequations}
This is essentially a problem for the Laplacian
with a generalized Robin boundary condition on~$\hat{\Gamma}$ involving a surface Laplacian (Laplace–Beltrami operator).  Again, to guarantee existence of a unique weak solution to this boundary-value problem, we use the Lax-Milgram theorem and establish coercivity of a suitable bilinear form. We assume~$\hat{\Omega}$ is a bounded Lipschitz-continuous domain, and $\hat{\phi}$ is sufficiently smooth to ensure continuity of the bilinear form and linear form of the weak formulation.

A weak form of~\eqref{DCS_L}--\eqref{DCS_R} seeks $\delta\phi \in V:=\{v\in H^1(\hat{\Omega})\,\big|\,v|_{\hat{\Gamma}}\in H^1(\hat{\Gamma})\}$ such that
\begin{equation}
    a(\delta\phi, \psi) = l(\psi), \quad \forall \psi\in V,
\end{equation}
where 
\begin{alignat}{2}
\label{a_uv}
a(\delta\phi,\psi) &:= \int_{\hat{\Omega}}\nabla\delta\phi\cdot\nabla \psi d\Omega + \int_{\hat{\Omega}\setminus\hat{\Gamma}}\omega\delta\phi \psi d\Gamma \\
&\qquad+\int_{\hat{\Gamma}}(\nabla_{\Gamma}\delta\phi)^T\big[\boldsymbol{b}_{\Gamma}(\nabla_{\Gamma}\hat{\phi})^T\big] \nabla_{\Gamma}\psi d\Gamma - \int_{\hat{\Gamma}}f\boldsymbol{b}_{\Gamma}\cdot\nabla_{\Gamma}\delta\phi \psi d\Gamma,
\nonumber
    \\
    \label{linear_right}
     l(v) &:= -\int_{\hat{\Omega}}\nabla\hat{\phi}\cdot\nabla \psi d\Omega + \int_{\hat{\Omega}}f\psi d\Omega \\
     &\qquad+ \int_{\partial\hat{\Omega}\setminus\hat{\Gamma}}(g+\omega h-\omega\hat{\phi})\psi d\Gamma + \int_{\hat{\Gamma}}\left[\text{div}_{\Gamma}\left(C\nabla_{\Gamma}\hat{\phi}\right)+Cf\right]\psi d\Gamma.\nonumber
\end{alignat}
To study the coercivity of $a(\cdot,\cdot)$, note that  
\begin{eqnarray}
    a(\psi,\psi) &=& \int_{\Omega}|\nabla\psi|^2d\Omega+\int_{\hat{\Omega}\setminus\hat{\Gamma}}\omega\psi^2d\Gamma \nonumber\\
    &&+\int_{\hat{\Gamma}}(\nabla_{\Gamma}\psi)^T\big[\boldsymbol{b}_{\Gamma}(\nabla_{\Gamma}\hat{\phi})^T\big]\nabla_{\Gamma}\psi d\Gamma-\int_{\hat{\Gamma}}f\boldsymbol{b}_{\Gamma}\cdot\nabla_{\Gamma}(\psi)\psi d\Gamma,
\end{eqnarray}
hence, similar to the approach in Appendix \ref{app:solvability_dirichlet}, if the last three terms are suitably bounded, coercivity holds. We note that the last two terms can be written as 
\begin{alignat*}{2}
 \int_{\hat{\Gamma}}(\nabla_{\Gamma}\psi)^T\big[\boldsymbol{b}_{\Gamma}(\nabla_{\Gamma}\hat{\phi})^T\big]\nabla_{\Gamma}\psi d\Gamma &= \int_{\hat{\Gamma}}(\nabla_{\Gamma}\psi)^T\big[\frac{2a\nabla_{\Gamma}\hat{\phi}
 }{C_1}(\nabla_{\Gamma}\hat{\phi})^T\big]\nabla_{\Gamma}\psi d\Gamma \\
 &= \int_{\hat{\Gamma}}\frac{2a
 }{C_1}(\nabla_{\Gamma}\hat{\phi}\cdot\nabla_{\Gamma}\psi)^2 d\Gamma 
    \\
\int_{\hat{\Gamma}}f\boldsymbol{b}_{\Gamma}\cdot\nabla_{\Gamma}(\psi)\psi d\Gamma &=\int_{\hat{\Gamma}}\frac{1}{2}f\boldsymbol{b}_{\Gamma}\cdot\nabla_{\Gamma}(\psi^2) d\Gamma 
\\ &= \int_{\hat{\Gamma}}\text{div}_{\Gamma}\left(\frac{1}{2}f\boldsymbol{b}_{\Gamma}\right)\psi^2d\Gamma - \int_{\partial\hat{\Gamma}}\psi^2(\frac{1}{2}f\boldsymbol{b}_{\Gamma}\cdot\boldsymbol{\tau})ds.
\end{alignat*}
Therefore, sufficient conditions that guarantee coercivity are:
\begin{equation}
\begin{cases}
    \frac{a}{C_1}\geq a_0>0, \quad \text{div}_{\Gamma}\left(f\boldsymbol{b}_{\Gamma}\right)\leq0,\quad f\boldsymbol{b}_{\Gamma}\cdot\boldsymbol{\tau}\geq0 \quad  \text{on } \hat{\Gamma}\\
    \omega>0 \quad \text{on } \partial\hat{\Omega}\setminus\hat{\Gamma}.
\end{cases}
\end{equation}
\par
The same remark at the end of Appendix~\ref{app:solvability_dirichlet} applies: The above conditions guarantee existence of $\delta \phi \in V$ while $\delta \boldsymbol{\theta} \cdot \boldsymbol{n} $ follows from~\eqref{linBC}. Generally, one expects additional regularity for $\delta \phi$ beyond~$V$, so that $\delta \boldsymbol{\theta}$ inherits this regularity and becomes a Lipschitz-continuous vector field on~$\Gamma_0$. Such regularity study is outside the scope of this work.

\section{Solvability of the discrete shape-linearized system for the Bernoulli condition}\label{app:solvability_discrete}
In this Appendix we show that, under certain conditions of the data and mesh, the finite element method for the shape-linearized system~\eqref{eq:discreteWeakForm} has a unique discrete solution. To that end, we use the Lax-Milgram theorem and establish coercivity of the coupled system in~\eqref{eq:discreteWeakForm}.\footnote{We note that the proof of coercivity in the continuous setting, see Appendix~\ref{app:solvability_bernoulli}, does not apply in the discrete case, because in the continuous setting the geometrical variable $\delta\boldsymbol{\theta} \cdot \boldsymbol{n}$ could be straightforwardly eliminated from the system. We therefore establish coercivity separately in the discrete case.}
 \par
The discrete problem \eqref{eq:discreteWeakForm} can be written as follows:
\begin{equation}
\begin{cases}
    \text{Find } (\delta\eta_h,\delta\phi_h)\in \hat{W}_h \times \hat{V}_h = \mathbb{P}^1_{0,\mathrm{in}}(\Gamma_h) \times \mathbb{P}^1(\Omega_h) :\\
    \qquad 
    a\left((\delta\eta_h,\delta\phi_h),(v_h,w_h)\right) = l(w_h,v_h) 
    \qquad \forall (w_h,v_h) \in \hat{W}_h\times \hat{V}_h\,,
\end{cases}\label{discrete_problem}
\end{equation}
where
\begin{subequations}
\begin{alignat*}{2}
 a((\delta\eta_h,\delta\phi_h), (w_h,v_h))&:= b\left((\delta\eta_h,\delta\phi_h),v_h\right) + c\left((\delta\eta_h,\delta\phi_h),w_h\right)
\\
 l(w_h,v_h) &:= l_{1}(v_h) +l_{2}(w_h)
\\
b\left((\delta\eta_h,\delta\phi_h),v_h\right)
&:=\int_{\hat{\Omega}}\nabla\delta\phi_h\cdot\nabla v_h d \Omega + \int_{\Omega\setminus\hat{\Gamma}}\omega\delta\phi_h v_hd\Gamma- \int_{\hat{\Gamma}}fv_h\delta\eta_h d \Gamma\nonumber\\
&\quad + \int_{\hat{\Gamma}}\nabla_{\Gamma
    }\hat{\phi}\cdot\nabla_{\Gamma}v_h \delta\eta _hd\Gamma ,\\
    l_{1}(v_h) &:=  {}-\mathcal{R}_1\left(\left(\hat{\boldsymbol\theta}, \hat{\phi}\right); v_h\right),\\
    c\left((\delta\eta_h,\delta\phi_h),w_h\right) &:=\int_{\hat{\Gamma}}2a\nabla_{\Gamma}\hat{\phi}\cdot\nabla_{\Gamma}\delta\phi_h w_h d \Gamma + \int_{\hat{\Gamma}}\left(2\kappa\left|\nabla\hat{\phi}\right|^2+bn_N\right)w_h\delta\eta_h d \Gamma, \\
    l_{2}(w_h)&:= {}-\mathcal{R}_2\left(\left(\hat{\boldsymbol\theta}, \hat{\phi}\right); w_h\right).
\end{alignat*}
\end{subequations}
To study the coercivity of $a(\cdot,\cdot)$, note that
\begin{subequations}
\begin{alignat}{2}
\label{ah}b\left((\delta\eta_h,\delta\phi_h),\delta\phi_h\right)=&\int_{\hat{\Omega}}\left|\nabla\delta\phi_h\right|^2 d \Omega + \int_{\partial\hat{\Omega}\setminus\hat{\Gamma}}\omega(\delta\phi_h)^2 d\Gamma \\
&+ \int_{\hat{\Gamma}}\nabla_{\Gamma
    }\hat{\phi}\cdot\nabla_{\Gamma}\delta\phi_h \delta\eta_hd\Gamma - \int_{\hat{\Gamma}}f\delta\phi_h\delta\eta_h d \Gamma,\nonumber\\
    c\left((\delta\eta_h,\delta\phi_h),\delta\eta_h\right) =&\int_{\hat{\Gamma}}2a\nabla_{\Gamma}\hat{\phi}\cdot\nabla_{\Gamma}\delta\phi_h \delta\eta_h d \Gamma + \int_{\hat{\Gamma}}\left(2\kappa\left|\nabla\hat{\phi}\right|^2+bn_N\right)(\delta\eta_h)^2 d \Gamma.\label{bh}
\end{alignat}
\end{subequations}

Next, we bound the coupling terms in~\eqref{ah} and~\eqref{bh}. Let $\hat{c}_0:= ||\nabla_{\Gamma}\hat{\phi}||_{L^\infty(\hat{\Gamma})}$ and $c_f:=||f||_{L^\infty(\hat{\Gamma})}$, then
\begin{alignat*}{2}
\left|\int_{\hat{\Gamma}}\nabla_{\Gamma}\hat{\phi}\cdot\nabla_{\Gamma}\delta\phi_{h}\delta\eta_h d\Gamma\right|
&\geq-\hat{c}_0||\nabla_{\Gamma}\delta\phi_h||_{L^2(\hat{\Gamma})}||\delta\eta_h||_{L^2(\hat{\Gamma})}\,,
\\
\left| \int_{\hat{\Gamma}}f\delta\phi_h\delta\eta_h d \Gamma \right |
&\geq  - c_f ||\delta\phi_h||_{L^2(\hat{\Gamma})}
||\delta\eta_h||_{L^2(\hat{\Gamma})}
\,.
\end{alignat*}
According to discrete trace inequalities~\cite[Ch. 12.2]{ern2021finite1}, there are constants~$c_1,c_2>0$ (independent of~$h$) such that
\begin{subequations}
\begin{alignat*}{2}
    ||\nabla_{\Gamma}\delta\phi_h||_{L^2(\hat{\Gamma})}
    \leq 
    ||\nabla \delta\phi_h||_{L^2(\hat{\Gamma})}
    &\leq
    c_1 h^{-\frac{1}{2}}||\nabla \delta\phi_h||_{L^2(\hat{\Omega})},
\\
||\delta\phi_h||_{L^2(\hat{\Gamma})}
    &\leq
    c_2 h^{-\frac{1}{2}}||\delta\phi_h||_{L^2(\hat{\Omega})},
\end{alignat*}
\end{subequations}
where $h = \max_{K\in \Omega_h} \mathrm{diam}(K) $.

Using these inequalities as well as Young's inequality (see, e.g.,~\cite[Appendix C.3]{ern2021finite}), we obtain
\begin{alignat}{2}
\left|\int_{\hat{\Gamma}}\nabla_{\Gamma}\hat{\phi}\cdot\nabla_{\Gamma}\delta\phi_{h}\delta\eta_h
d\Gamma\right|
&\geq -\frac{\hat{c}_0 c_1^2 h^{-1}}{2} ||\nabla\delta\phi_h||_{L^2(\hat{\Omega})}^2- \frac{\hat{c}_0}{2}||\delta\eta_h||_{L^2(\hat{\Gamma})}^2,\label{appd:inequality_1}
\\
\left| \int_{\hat{\Gamma}}f\delta\phi_h\delta\eta_h d \Gamma \right |
&\geq  - \frac{c_f c_2^2 h^{-1}}{2} ||\delta\phi_h||_{L^2(\hat{\Omega})}^2
- \frac{c_f }{2}
||\delta\eta_h||^2_{L^2(\hat{\Gamma})}
\,.\label{appd:inequality_2}
\end{alignat}
Next, using Poincare inequality,
\begin{equation}
||\delta\phi_h||^2_{L^2(\hat{\Omega})}\leq c_p\left(||\nabla\delta\phi_h||^2_{L^2(\hat{\Omega})}+||\delta\phi_h||^2_{L^2(\partial\hat{\Omega}\setminus\hat{\Gamma})}\right)\label{appd:poincare}
\end{equation}
And substitution from \eqref{appd:inequality_1}-\eqref{appd:poincare} into \eqref{ah} and \eqref{bh}, we have
\begin{subequations}
\begin{alignat}{2}
b\left((\delta\eta_h,\delta\phi_h),\delta\phi_h\right)\geq&\left(1-\frac{\hat{c}_0 c_1^2 h^{-1}}{2}-\frac{c_pc_fc_2^2h^{-1}}{2}\right)||\nabla\delta\phi_h||_{L^2(\hat{\Omega})}^2 \\
&+ \left(\omega-\frac{c_pc_fc_2^2h^{-1}}{2}\right)||\delta\phi_h||^2_{L^2(\partial\hat{\Omega}\setminus\hat{\Gamma})}
-\left(\frac{\hat{c}_0+c_f}{2}\right)||\delta\eta_h||^2_{L^2(\hat{\Gamma})}\nonumber\\
c\left((\delta\eta_h,\delta\phi_h),\delta\eta_h\right) =&-|a|\hat{c}_0 c_1^2 h^{-1} ||\nabla\delta\phi_h||_{L^2(\hat{\Omega})}^2\\
&+\left(\left|2\kappa|\nabla\hat{\phi}|^2+bn_{N}\right|- |a|\hat{c}_0\right)||\delta\eta_h||_{L^2(\hat{\Gamma})}^2.\nonumber
\end{alignat}
\end{subequations}
\par
Hence, we obtain the estimate
\begin{alignat}{2}
a\left((\delta\eta_h,\delta\phi_h),(\delta\eta_h,\delta\phi_h)\right) &\geq 
    \left(1-\frac{\hat{c}_0c_1^2 h^{-1} |1+2a|}{2}-\frac{c_pc_fc_2^2h^{-1}}{2}\right)||\nabla\delta\phi_h||^2_{L^2(\hat{\Omega})}\nonumber\\
    &+\left(\omega-\frac{c_pc_fc_2^2h^{-1}}{2}\right)||\delta\phi_h||^2_{L^2(\partial\hat{\Omega}\setminus\hat{\Gamma})}\nonumber\\
    &+\left(\left|2\kappa\nabla\hat{\phi}|^2+bn_N\right|-\frac{\left(|1+2a|\hat{c}_0+c_f\right)}{2}\right)||\delta\eta_h||^2_{L^2(\hat{\Gamma})},
\end{alignat}
Therefore, sufficient conditions that guarantee discrete solvability are:
\begin{subequations}
    \begin{gather*}
        1-\frac{\hat{c}_0c_1^2 h^{-1} |1+2a|}{2}-\frac{c_pc_fc_2^2h^{-1}}{2}>0, \\
        \omega-\frac{c_pc_fc_2^2h^{-1}}{2}>0, \\
\left|2\kappa\nabla\hat{\phi}|^2+bn_N\right|-\frac{\left(|1+2a|\hat{c}_0+c_f\right)}{2}>0.
    \end{gather*}
\end{subequations}

\section*{Acknowledgments}
The authors are grateful to Anna Kalogirou and Onno Bokhove for additional discussion. The authors would also like to thank the anonymous reviewers for their helpful comments and suggestions, which led to many significant improvements, in particular the addition of Section~\ref{section_3d} (on 3-D equivalent expression) and Remarks~\ref{remark_solve_dirichlet} and~\ref{remark_solve_discrete}, and corresponding Appendices (on solvability of the continuous and discrete linearised systems).

\bibliographystyle{siamplain}
\bibliography{references}
\end{document}